\documentclass[a4paper,10pt]{article}
\usepackage[top=13mm, bottom=23mm, inner=20mm, outer=20mm]{geometry}
\title{A well-balanced reconstruction with bounded velocities for the shallow water equations by convex combination}
\author{Edward W. G. Skevington}
\date{21/06/2021}

\usepackage{mathtools,amssymb,mathrsfs,bm}

\usepackage{graphicx}
\usepackage{caption,subcaption,float}
\usepackage{tabularx,multirow}
\usepackage[inline]{enumitem}

\usepackage{xcolor}

\numberwithin{equation}{section}
\numberwithin{figure}{section}
\numberwithin{table}{section}

\usepackage{calc,etoolbox}

\usepackage[inline]{enumitem}

\makeatletter
\@ifpackageloaded{hyperref}{}{\usepackage[hidelinks]{hyperref}}
\@ifpackageloaded{cleveref}{}{
	\usepackage[compress]{cleveref}

	\crefformat{section}{\S#2#1#3}
	\crefformat{subsection}{\S#2#1#3}
	\crefformat{subsubsection}{\S#2#1#3}
	\crefrangeformat{section}{\S\S#3#1#4 to~#5#2#6}
	\crefmultiformat{section}{\S\S#2#1#3}{ and~#2#1#3}{, #2#1#3}{ and~#2#1#3}
	
	\crefname{figure}{fig.}{figs.}
	\Crefname{figure}{Fig.}{Figs.}
	
	\crefname{equation}{}{}
	\crefformat{equation}{(#2#1#3)}
	\crefmultiformat{equation}{(#2#1#3)}{ and (#2#1#3)}{, (#2#1#3)}{ and (#2#1#3)}
}
\makeatother

\usepackage[title]{appendix}
\newenvironment{myappendices}
{% \begin{myappendices}
	\begin{appendices}
	\crefalias{section}{appsec}
}{% \end{myappendices}
	\end{appendices}
}

\crefname{appsec}{Appendix}{Appendices}
\Crefname{appsec}{Appendix}{Appendices}

\usepackage{xspace}
\makeatletter
\DeclareRobustCommand\onedot{\futurelet\@let@token\@onedot}
\def\@onedot{\ifx\@let@token.\else.\null\fi\xspace}
\makeatother
\def\eg{e.g\onedot} 
\def\ie{i.e\onedot} 
\def\cf{\textit{cf}\onedot}

\def\etal{et al\onedot}
\def\sf{s.f\onedot}

\makeatletter
\newcommand{\floattag}[1]{%
  \@namedef{the\@captype}{#1}%
  \@namedef{theH\@captype}{#1}%
  \addtocounter{\@captype}{-1}}
\makeatother
\input{Templates/StandardMathsMacros}
\usepackage{amsthm, thmtools, thm-restate}
\theoremstyle{plain}
\declaretheorem[numberwithin=section,	refname={theorem,theorems},			Refname={Theorem,Thorems}]								{theorem}
\declaretheorem[name=Lemma,				refname={lemma,lemmas},				Refname={Lemma,Lemmas},				sibling=theorem]	{lemma}

\theoremstyle{definition}
\declaretheorem[name=Definition,		refname={definition,definitions},	Refname={Definition,Definitions},	sibling=theorem]	{definition}

\theoremstyle{remark}

\theoremstyle{definition}

\crefformat{testproblem}{test problem #2#1#3}
\crefrangeformat{testproblem}{test problems #3#1#4 to~#5#2#6}
\Crefrangeformat{testproblem}{Test problems #3#1#4 to~#5#2#6}
\crefmultiformat{testproblem}{test problems #2#1#3}{ and~#2#1#3}{, #2#1#3}{ and~#2#1#3}

\graphicspath{{Figures/}}

\allowdisplaybreaks[1]

\usepackage{printlen}
\begin{document}

\maketitle

\begin{abstract}
Finite volume schemes for hyperbolic balance laws require a piecewise polynomial reconstruction of the cell averaged values, and a reconstruction is termed `well-balanced' if it is able to simulate steady states at higher order than time evolving states. For the shallow water system this involves reconstructing in surface elevation, to which modifications must be made as the fluid depth becomes small to ensure positivity, and for many reconstruction schemes a modification of the inertial field is also required to ensure the velocities are bounded. 

We propose here a reconstruction based on a convex combination of surface and depth reconstructions which ensures that the depth increases with the cell average depth. We also discuss how, for cells that are much shallower than their neighbours, reducing the variation in the reconstructed flux yields bounds on the velocities. This approach is generalisable to high order schemes, problems in multiple spacial dimensions, and to more complicated systems of equations. We present reconstructions and associated technical results for three systems, the standard shallow water equations, shallow water in a channel of varying width, and a shallow water model of a particle driven current. Positivity preserving time stepping is also discussed.
\end{abstract}

%\begin{keywords}
%Finite Volume, Conservation Laws, Shallow Water, Well-Balanced
%\end{keywords}
%
%\begin{AMS}
%65M08, 76M12, 35L65
%\end{AMS}
\section{Introduction}	\label{sec:Intro}

The construction of \emph{well-balanced} numerical schemes has been the subject of extensive research for the shallow water equations (\eg \cite{bk_Stoker_WW,bk_Ungarish_GCI}),
\begin{subequations}	\label{eqn:INTRO_SW_Full}	\begin{align}
	\pdv{h}{t} + \pdv{}{x}(uh) &= 0,		\label{eqn:INTRO_SW_Full_h}	\\
	\pdv{}{t}(uh) + \pdv{}{x}\ppar*{u^2 h + \frac{g h^2}{2}} &= - g h\dv{b}{x},			\label{eqn:INTRO_SW_Full_u}
\end{align}\end{subequations}
where $h(x,t)$ is the fluid depth, $b(x)$ the bed elevation, $u(x,t)$ the horizontal velocity, and $g$ the gravitational acceleration. The system may be modified to include additional effects, for example including source term in \cref{eqn:INTRO_SW_Full_h} to model fluid entrainment, or in \cref{eqn:INTRO_SW_Full_u} to capture viscous effects, or the inclusion a density field (\eg \cite{bk_Ungarish_GCI}). 
%In addition, \cref{eqn:INTRO_SW_Full} is only in one spatial dimension (not counting depth), and many natural phenomena occur in two dimensions.

\begin{figure}[t]
	\begin{center}
%% Accompanies image file 'DrainIntoLake.pdf' (pdf, eps, ps)
\begingroup%
  \makeatletter%
  \providecommand\color[2][]{%
    \errmessage{(Inkscape) Color is used for the text in Inkscape, but the package 'color.sty' is not loaded}%
    \renewcommand\color[2][]{}%
  }%
  \providecommand\transparent[1]{%
    \errmessage{(Inkscape) Transparency is used (non-zero) for the text in Inkscape, but the package 'transparent.sty' is not loaded}%
    \renewcommand\transparent[1]{}%
  }%
  \providecommand\rotatebox[2]{#2}%
  \newcommand*\fsize{\dimexpr\f@size pt\relax}%
  \newcommand*\lineheight[1]{\fontsize{\fsize}{#1\fsize}\selectfont}%
  \ifx\svgwidth\undefined%
    \setlength{\unitlength}{325.98425197bp}%
    \ifx\svgscale\undefined%
      \relax%
    \else%
      \setlength{\unitlength}{\unitlength * \real{\svgscale}}%
    \fi%
  \else%
    \setlength{\unitlength}{\svgwidth}%
  \fi%
  \global\let\svgwidth\undefined%
  \global\let\svgscale\undefined%
  \makeatother%
  \begin{picture}(1,0.26086957)%
    \lineheight{1}%
    \setlength\tabcolsep{0pt}%
    \put(0,0){\includegraphics[width=\unitlength,page=1]{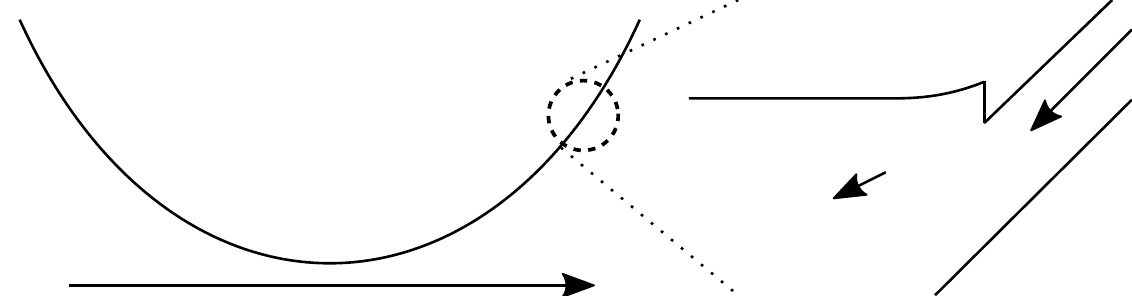}}%
    \put(0.51024246,0.01784084){\color[rgb]{0,0,0}\makebox(0,0)[lt]{\lineheight{1.25}\smash{\begin{tabular}[t]{l}$x$\end{tabular}}}}%
    \put(0.38542616,0.05262345){\color[rgb]{0,0,0}\makebox(0,0)[rt]{\lineheight{1.25}\smash{\begin{tabular}[t]{r}$b$\end{tabular}}}}%
    \put(0.28741057,0.18305797){\color[rgb]{0,0,0}\makebox(0,0)[t]{\lineheight{1.25}\smash{\begin{tabular}[t]{c}$\eta$\end{tabular}}}}%
    \put(0.76920047,0.11101976){\color[rgb]{0,0,0}\makebox(0,0)[rt]{\lineheight{1.25}\smash{\begin{tabular}[t]{r}$u$\end{tabular}}}}%
    \put(0,0){\includegraphics[width=\unitlength,page=2]{DrainIntoLake.pdf}}%
  \end{picture}%
\endgroup%

	\end{center}
	\caption{Illustration of a typical situation, where the flow is approaching steady state leaving some fluid still on the slopes. The body of the fluid is close to `lake at rest', and on the slopes the fluid is close to `thin film'.}
	\label{fig:PWB_draining_flow}
\end{figure}

A numerical scheme simulating a system of hyperbolic balance laws is termed \emph{well-balanced} \cite{ar_Greenberg_1996} (also known as the C-property \cite{ar_Bermudez_1994}) if it resolves steady states at higher accuracy than time evolving states, and in some cases may be accurate to machine precision. This allows simulations to capture solutions that are perturbations of, or tend towards, a steady state. For \cref{eqn:INTRO_SW_Full}, a special steady state is $h = 0$, the \emph{dry bed}. All other steady states are of constant volume flux, $q$, and energy, $E$, where
\begin{align}	\label{eqn:SW_steady_moving}
	q &\eqdef uh,
	&
	E &\eqdef \frac{u^2}{2} + g \eta,
\end{align}
and $\eta \eqdef h+b$ is the surface elevation. This state is termed the \emph{moving water equilibrium}, and it has two important special cases illustrated in \cref{fig:PWB_draining_flow}. Firstly the \emph{lake at rest} state with $u=0$ and $\eta = \bar{\eta}$ a constant, which exists for a wide class of generalized shallow water systems. Secondly the \emph{thin film} state where the fluid is draining off a slope with $\dv*{b}{x} = \Torder{\frac*{H}{L}}$ so that $h \ll H$. Here $H$ and $L$ are the vertical and horizontal length-scales respectively, and we employ asymptotic notation $\Torder{}$. In this state $u = \Torder{\sqrt{gH}}$, thus $\abs{\frac*{u}{\sqrt{gh}}} = \Torder{\sqrt{\frac*{H}{h}}} \gg 1$ and
\begin{equation} \label{eqn:thin_film}
	\dv{h}{x} = \frac{1}{(u^2/gh) - 1} \dv{b}{x} = \Torder{\frac{h}{L}}.
\end{equation}
The variation in depth is small, while the variation in surface elevation is large, $\dv*{\eta}{x} = \Torder{H/L}$. 
We expect the thin film regime to result in slowly varying $h$ for a wide range of shallow water models, in particular including basal friction makes the correct steady state constant $h$ \cite{ar_Chertock_2015}. 
If we wish to converge rapidly on exact solutions then we need to resolve thin film regions accurately. For example, a dam break with a tailwater of thickness $h_t \ll 1$ (\eg \cite{ar_Delestre_2016}) converges on the $h_t = 0$ solution at $\Torder{h_t^{1/4}}$. We expect that the propagation of fluid up a slope \cite{ar_Carrier_1958,ar_Thacker_1981} will converge similarly slowly with respect to thin films, thus the dynamics of the thin film itself should be resolved.

We present a reconstruction for finite volume schemes that operates in $\eta$ where the fluid is deep (relative to bed variation) and in $h$ where it is shallow or fast. Our reconstruction differs from many present in the literature \cite{ar_Kurganov_2002,ar_Bollermann_2013,ar_Chertock_2015} in that, as the fluid becomes deeper, the reconstructed depths also become deeper, a property we term self-monotonicity (\cref{sec:overview}). We discuss three different systems: classic shallow water \cref{eqn:INTRO_SW_Full}; channel flow with varying width \cref{eqn:SW_width_Full}; and fluid driven by a settling particle load \cref{eqn:EXT_SW_Full}, for which we also discuss positivity preserving time stepping (\cref{sec:postimestep}). Our approach is readily generalisable to other systems and higher dimensions.

The paper is structured as follows. We begin with a brief summary of finite volume schemes in \cref{sec:FiniteVolume}, before discussing our approach to reconstruction in \cref{sec:overview}. In \cref{sec:results} we present a reconstruction which uses our approach and the minmod slope limiter, including bounds on the reconstructed depths and velocities, which is contrasted with other reconstructions in \cref{sec:comparison}. In \cref{sec:monotone} self-monotonicity is discussed for scalar problems, We finish by presenting proofs for our theorems in order throughout \cref{sec:reconH,sec:recon_q,sec:reconWH,sec:recon_phi}, one theorem per section. As an associated result, we include a method of positivity preserving time evolution in \cref{sec:postimestep}. The supplemental material proves more general results than those detailed in \cref{sec:results}, outlining of how similar results may be proven using other slope limiters or higher order alternatives.
\section{Finite volume schemes} \label{sec:FiniteVolume}

%A variety of numerical approaches exist for hyperbolic balance laws, \eg well-balanced discontinuous Galerkin schemes \cite{ar_Noelle_2009,ar_Xing_2010,ar_Bollermann_2011,ar_Xing_2014,ar_Chen_2019}. We consider finite volume schemes (\eg \cite{bk_Leveque_FVM}), which we briefly introduce for the general balance law

Consider the balance law
\begin{equation}	\label{eqn:FV_conservation}
	\pdv{Q}{t} + \pdv{}{x} (F) = \Psi	\quad \textrm{on} \quad x_L \leq x \leq x_R,
\end{equation}
where $Q$ is a function from $(x,t)$ to $\mathbb{R}^M$, and the flux $F$ and source $\Psi$ are functions from $(Q,x,t)$ to $\mathbb{R}^M$. The system is spatially discretized over $J$ cells by introducing $J+1$ cell interface points $x_L=x_{1/2}<x_{3/2}<\ldots<x_{J+1/2}=x_R$ where the $j^\nth$ cell $C_j = (x_{j-1/2},x_{j+1/2})$ has width $\Delta x_j=x_{j+1/2} - x_{j-1/2}$, and
\begin{align}\label{eqn:FV_QPdef} 
	Q_{j}(t) &\eqdef \frac{1}{\Delta x_j} \int_{C_j} Q(x,t) \dd{x},
	&\textrm{and}&&
	\Psi_{j}(t) &\eqdef \frac{1}{\Delta x_j} \int_{C_j} \Psi(Q(x,t),x,t) \dd{x}
\end{align}
are cell averaged values.  Averaging \cref{eqn:FV_conservation} over $C_j$ yields
\begin{equation}
	\dv{Q_{j}}{t} + \frac{1}{\Delta x_j} \left[ F_{j+1/2} - F_{j-1/2} \right] = \Psi_{j},
\end{equation}
where $F_{j+1/2}(t) \eqdef F(Q(x_{j+1/2},t),x_{j+1/2},t)$. A finite volume scheme approximates $F_{j+1/2}$ and $\Psi_j$, and often this in done in two steps \cite{ar_vanLeer_1979,ar_Osher_1985}. First, $Q$ is \emph{reconstructed} by a polynomial over each cell, yielding $\hat{Q}(x,t)$ with discontinuities at the cell interfaces
\begin{align}\label{eqn:FV_recon} 
	Q_{j-1/2}^+(t) &\eqdef \lim_{x \rightarrow x_{j-1/2}^+} \hat{Q}(x,t),
	&
	Q_{j+1/2}^-(t) &\eqdef \lim_{x \rightarrow x_{j+1/2}^-} \hat{Q}(x,t).
\end{align}
Second, the numerical flux $F_{j+1/2} = \tilde{F}(Q_{j+1/2}^-,Q_{j+1/2}^+,x_{j+1/2},t)$ is deduced: in Godunov schemes, a Riemann problem at the cell interface is solved (either exactly or approximately) to find the solution at $(x_{j+1/2},t^+)$, denoted $Q_{j+1/2}^R$, and then $F_{j+1/2} =  F(Q_{j+1/2}^R,x_{j+1/2},t)$; while in central schemes the flux is computed directly from the reconstruction. The choice of reconstruction, flux, source, and time stepping algorithm determines the scheme.
\section{Self-monotone well-balancing by convex combination} \label{sec:overview}

A number of well-balanced schemes have previously been developed. Some rely on specially designed expression for the flux \cite{ar_Berthon_2008,ar_Berthon_2016}, while others reconstruct the surface elevation and depth and deduce the bed structure \cite{ar_Audusse_2004,ar_Noelle_2006}, or use somewhat exotic properties of the specific system \cite{ar_Audusse_2005}. Of interest here are those that are well-balanced by the choice of reconstruction and the approximation of the source term. One approach is to reconstruct in flux and energy \cref{eqn:SW_steady_moving} \cite{ar_Xing_2014,ar_Cheng_2016}, which is challenging to generalise to other similar systems, such as the case of shallow water with basal friction \cite{ar_Cheng_2019} or higher dimensions. Others focus on the lake at rest state. If the discretized $\eta$ is constant and $q$ is zero, \ie $\eta_j = \bar{\eta}$, $q_j = 0$, then it is reasonable that the reconstruction should also be constant, \ie $\eta_{j+1/2}^\pm = \bar{\eta}$, $q_{j+1/2}^\pm=0$. Then, using the discrete expression for the source term from \cite{ar_Bermudez_1994} which assumes $b$ continuous (used in \cite{ar_Jin_2001,ar_Kurganov_2002,ar_Audusse_2004,ar_Noelle_2006,ar_Kurganov_2007,ar_Bollermann_2011}, generalised to a discontinuous bed in \cite{ar_Bernstein_2016})
\begin{equation}	\label{eqn:source_discretisation}
	-gh \dv{b}{x} 
	\approx
	-g \frac{h_{j-1/2}^+ + h_{j+1/2}^-}{2} \frac{b_{j+1/2} - b_{j-1/2}}{\Delta x_j}
\end{equation}
(for Godunov schemes $h_{j+1/2}^\pm$ may be replaced with $h_{j+1/2}^R$) a balance between fluxes and source terms is obtained
\begin{align*}
	\frac{F_{j+1/2} - F_{j-1/2}}{\Delta x_j}
	=
	\frac{g}{2 \Delta x_j}
	\pgrp*{\begin{array}{c}
		0 \\
		(b_{j+1/2})^2 - 2 \bar{\eta} b_{j+1/2} + 2 \bar{\eta} b_{j-1/2} - (b_{j-1/2})^2
	\end{array}}
	=
	\Psi_j.
\end{align*}
By construction, \cref{eqn:source_discretisation} limits the scheme to second order accuracy \cite{ar_Kurganov_2002}. For consistency, schemes employing \cref{eqn:source_discretisation} typically use a piecewise linear reconstruction.

The main difference between the various schemes is how they enforce the \emph{positivity preserving} property, $h_j \geq 0$ and $h_{j+1/2}^\pm \geq 0$, which is required for the solution to be physical and the system to be hyperbolic. This is only a concern when the variation in $h$ or $\eta$ becomes of order the variation of $b$, such as close to a transition from `lake at rest' to `dry bed' or `thin film' (\cref{fig:PWB_draining_flow}). The earliest approach \cite{ar_Kurganov_2002} used a piecewise linear reconstruction in $\eta$ when the fluid exceeded a depth threshold in the cell and its neighbours, otherwise reconstructing in $h$. In \cite{ar_Kurganov_2007,ar_Balbas_2014,ar_Bernstein_2016} the same reconstruction was used in deep areas, but when the flow became shallow the only change was to enforce that, whenever the depth at a cell interface was negative (\eg $h_{j+1/2}^- < 0$), the reconstructed gradient would be adjusted so that the depth there was zero instead ($h_{j+1/2}^- = 0$). Setting the depth to precisely vanish caused problems when evaluating the velocity field, so an additional modification was made to the computation of $u$ by using $q/(h + \epsilon(h))$. In \cite{ar_Bollermann_2013}, any cell that contained insufficient fluid for a constant surface elevation was given a reconstruction consisting of two pieces, one bringing the depth down to zero, the other of constant zero depth. In \cite{ar_Chertock_2015} the reconstruction was modified in cells that contain negative reconstructed depths by setting $h$ to be constant within the cell.

To critique these approaches we define some terms for a scalar field $v$ (\eg $v=h$) with cell averaged values $v_j$. Reconstructing $\hat{v}$ yields cell interface values $v_{j+1/2}^{\pm}$ which are functions of the cell averages $\{v\} \eqdef (v_1 , v_2 \ldots v_J)$ and cell interfaces  $\{x\} \eqdef (x_{1/2} , x_{3/2} \ldots x_{J+1/2} )$ (and potentially the values of other fields also).
\begin{definition} \label{def:monotone}
	A reconstruction $\hat{v}$ is \emph{self-monotone} when it is a non-decreasing function of $v_j$ for $x \in C_j$. A reconstruction is \emph{neighbour-monotone} when $v_{j-1/2}^{-}$ and $v_{j+1/2}^{+}$ are non-decreasing functions of $v_{j}$. 
\end{definition}
\begin{lemma}
	A reconstruction that is linear over each cell is self-monotone precisely when $v_{j-1/2}^{+}$ and $v_{j+1/2}^{-}$ are non-decreasing functions of $v_j$.
\end{lemma}
The desirability of the properties we hold to be self evident, though for scalar problems it is possible to give a formal justification which we present in \cref{sec:monotone}.

The reconstructions from \cite{ar_Kurganov_2002,ar_Bollermann_2013,ar_Chertock_2015} are not self-monotone in depth (or, equivalently, in surface elevation). That is, if we take the discrete values of depth $\{h\} \eqdef (h_1 , h_2 \ldots h_J)$, and hold all but one constant, the remaining one $h_j$ we increase, then $\hat{h}$ may decrease in cell $C_j$. 
%Even high order, discretized Galerkin schemes may feature a non-self-monotone reconstruction \cite{ar_Xing_2010} \note{check}. 
The schemes from \cite{ar_Kurganov_2007,ar_Balbas_2014,ar_Bernstein_2016,ar_Bollermann_2013,ar_Chertock_2015} require modification of the velocity field, which is strange because they reconstruct in $q$ which is constant in steady state. This is partially due to a lack of a lower bound on depth.

Our new approach is to perform two reconstructions, one in $h$ denoted $\hat{h}^h$, the other in $\eta$, denoted $\hat{h}^\eta$, and then perform a \emph{convex combination} of the results to obtain
\begin{align}	\label{eqn:convex_combination}
	\hat{h} &\eqdef \gamma_j \hat{h}^{\eta} + (1-\gamma_j) \hat{h}^h 
	&&\textrm{for}&
	x &\in C_j
\end{align}
with $\gamma_j$ a function of $\{Q\} \eqdef (Q_1 , Q_2 \ldots Q_J)$, $\{b\} \eqdef (b_{1/2}^+ , b_{3/2}^- , b_{3/2}^+ \ldots b_{J+1/2}^-)$ and $\{x\} \eqdef (x_{1/2}, x_{3/2} \ldots x_{J+1/2})$ with codomain $[0,1]$. The convex combination allows us to smoothly transition from one mode of reconstruction to the other as the discrete depths change, rather than using a sharp cut-off \cite{ar_Kurganov_2002,ar_Xing_2010}, and thereby obtain a self-monotone reconstruction with a lower bound for the depth. Our reconstruction weakly violates neighbour-monotonicity, which we discuss further in \cref{sec:comparison,sec:supp_reconH}, though this is not believed to be a general limitation on \cref{eqn:convex_combination}.

\begin{figure}[tp!]
	%\centering
	\begin{center}
%% Accompanies image file 'UnboundedVelocity.pdf' (pdf, eps, ps)
\begingroup%
  \makeatletter%
  \providecommand\color[2][]{%
    \errmessage{(Inkscape) Color is used for the text in Inkscape, but the package 'color.sty' is not loaded}%
    \renewcommand\color[2][]{}%
  }%
  \providecommand\transparent[1]{%
    \errmessage{(Inkscape) Transparency is used (non-zero) for the text in Inkscape, but the package 'transparent.sty' is not loaded}%
    \renewcommand\transparent[1]{}%
  }%
  \providecommand\rotatebox[2]{#2}%
  \newcommand*\fsize{\dimexpr\f@size pt\relax}%
  \newcommand*\lineheight[1]{\fontsize{\fsize}{#1\fsize}\selectfont}%
  \ifx\svgwidth\undefined%
    \setlength{\unitlength}{198.42519685bp}%
    \ifx\svgscale\undefined%
      \relax%
    \else%
      \setlength{\unitlength}{\unitlength * \real{\svgscale}}%
    \fi%
  \else%
    \setlength{\unitlength}{\svgwidth}%
  \fi%
  \global\let\svgwidth\undefined%
  \global\let\svgscale\undefined%
  \makeatother%
  \begin{picture}(1,0.42857143)%
    \lineheight{1}%
    \setlength\tabcolsep{0pt}%
    \put(0,0){\includegraphics[width=\unitlength,page=1]{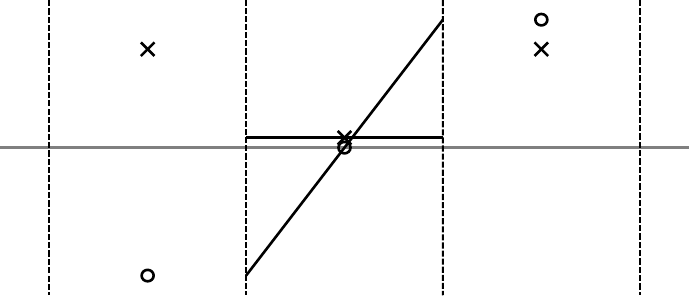}}%
    \put(0.61428571,0.38571429){\color[rgb]{0,0,0}\makebox(0,0)[rt]{\lineheight{1.25}\smash{\begin{tabular}[t]{r}$q$\end{tabular}}}}%
    \put(0.37142857,0.24285714){\color[rgb]{0,0,0}\makebox(0,0)[lt]{\lineheight{1.25}\smash{\begin{tabular}[t]{l}$h$\end{tabular}}}}%
    \put(0.5,0.01428571){\color[rgb]{0,0,0}\makebox(0,0)[t]{\lineheight{1.25}\smash{\begin{tabular}[t]{c}$C_j$\end{tabular}}}}%
  \end{picture}%
\endgroup%

	\end{center}
	\caption{An example reconstruction in cell $C_j$, with discrete cell average depths marked with crosses, and discrete volume flux with open circles, the line $h=q=0$ is in grey. Decreasing $h_j$ causes the reconstructed velocities to increase without bound.}
	\label{fig:unbounded_velocity}
\end{figure}
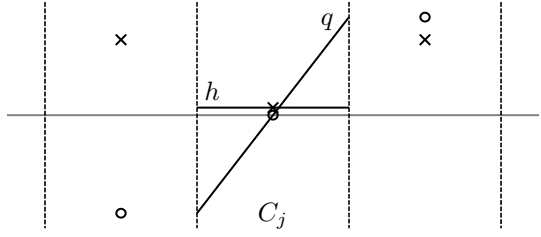

To obtain an upper bound for the velocities we first note that, even without bed variation, it is possible for the reconstructed velocities to become unbounded. We plot a worst case scenario in \cref{fig:unbounded_velocity}, where the depth in cell $C_j$ is vastly smaller than its neighbours and the sign of velocity changes sign across the cell. Then even if the reconstruction in depth is constant, as $h_j \rightarrow 0$ we will have an unbounded divergence of $u_{j-1/2}^+$ and $u_{j+1/2}^-$. To bound the velocities we \emph{suppress} the flux reconstruction as
\begin{align}	\label{eqn:flux_suppress}
	\hat{q} &= \kappa_j \hat{q}^q + (1-\kappa_j) q_j
	&&\textrm{for}&
	x &\in C_j
\end{align}
where $\hat{q}^q$ is a reconstruction using some standard method, and $\kappa_j$ is a function of $\{Q\}$ and $\{x\}$ with codomain $[0,1]$. We bound the velocities by choosing $\kappa_j$ to typically be $1$ but limits to $0$ when $h_j$ is much smaller than its neighbours. Convex combination and suppression are applicable to high order schemes and multiple spatial dimensions.
\section{Discussion of results} \label{sec:results}

To explore how a self-monotone reconstruction may be designed we focus on the minmod slope limiter employed by many well balanced schemes \cite{ar_Kurganov_2002,ar_Kurganov_2007,ar_Bollermann_2013,ar_Chertock_2015}, 
\begin{subequations}	\label{eqn:minmod_recon}
\begin{align}	\label{eqn:minmod_recon_R}	\begin{split}
	\sigma(v_{j-1},v_{j},v_{j+1}; \{x\}) &\eqdef \frac{2}{\Delta x_j}\minmod \left[
		\alpha_{j-1/2}^{+} (v_j - v_{j-1}),%\Delta v_{j-1/2} , 
	\right. \\	& \qquad \left.
		\alpha_{j} (v_{j+1}-v_{j-1}) , 
		\alpha_{j+1/2}^{-} (v_{j+1} - v_j)%\Delta v_{j+1/2}
	\right], 
\end{split}\end{align}
where $v$ is a general scalar field, $\alpha_j , \alpha_{j+1/2}^\pm : \{x\} \rightarrow [0,1]$ (\cref{thm:PWB_linearTVD}), and
\begin{align}\label{eqn:minmod_recon_MM}
	\minmod[y_1 , y_2 \ldots y_n] \eqdef
	\begin{cases}
		\min[y_1 , y_2 \ldots y_n] 	& \textrm{if } \min[y_1 , y_2 \ldots y_n]>0,	\\
		\max[y_1 , y_2 \ldots y_n] 	& \textrm{if } \max[y_1 , y_2 \ldots y_n]<0,	\\
		0							& \textrm{otherwise}.
	\end{cases}
\end{align}
\end{subequations}
The reconstruction of all fields is done using \cref{eqn:minmod_recon} with the same parameters. Any of the common expressions for $\minmod$ slope limiters may be obtained by selecting the parameters, \eg for the expression in \cite{ar_Kurganov_2019} take $\alpha_{j-1/2}^{+} = \frac*{\theta \Delta x_{j}}{( \Delta x_{j-1} + \Delta x_{j} )}$, $\alpha_{j-1/2}^{+} = \frac*{\Delta x_{j}}{( \Delta x_{j-1} + 2\Delta x_{j} + \Delta x_{j+1} )}$, $\alpha_{j+1/2}^{-} = \frac*{\theta \Delta x_{j}}{( \Delta x_{j} + \Delta x_{j+1} )}$.

Throughout this and the following sections we denote for all fields $v$
\begin{subequations}	\label{eqn:field_identities}
\begin{align}
	v_j &\equiv \frac{v_{j+1/2}^- + v_{j-1/2}^+}{2},
	&
	\Delta v_{j+1/2} &\equiv v_{j+1} - v_{j},
	\\
	\Delta v_j &\equiv v_{j+1/2}^- - v_{j-1/2}^+,
	&
	[v_x]_j &\equiv \frac{\Delta v_j}{\Delta x_j}.
\intertext{The functions $v_j^\downarrow$ and $v_j^\uparrow$ are computed as}
	v_j^\downarrow &\eqdef \min \pbrk*{ v_{j-1/2}^+ , v_{j+1/2}^- },
	&
	v_j^\uparrow   &\eqdef \max \pbrk*{ v_{j-1/2}^+ , v_{j+1/2}^- }
	\label{eqn:field_identities_vpm_known}
\end{align}
for fields where $v_{j+1/2}^\pm$ are known prior to reconstruction, \eg $b$ (note that, while $\alpha$ is not a field, we will use the notation \cref{eqn:field_identities_vpm_known}), and as
\begin{equation}\begin{aligned}
	v_j^\downarrow &\eqdef \min \pbrk*{ v_j - \alpha_{j-1/2}^{+}\Delta v_{j-1/2},v_j,v_j + \alpha_{j+1/2}^{-}\Delta v_{j+1/2} },
	\\
	v_j^\uparrow &\eqdef \max \pbrk*{ v_j - \alpha_{j-1/2}^{+}\Delta v_{j-1/2},v_j,v_j + \alpha_{j+1/2}^{-}\Delta v_{j+1/2} },
	\label{eqn:field_identities_vpm_recon}
\end{aligned}\end{equation}
for reconstructed fields, \ie $Q$, where $\alpha_{j+1/2}^\pm$ are parameters from \cref{eqn:minmod_recon}. Fields reconstructed using \cref{eqn:minmod_recon} satisfy
\begin{align}
	v_j^\downarrow &\leq \min \pbrk*{ v_{j-1/2}^+ , v_{j+1/2}^- },
	&
	v_j^\uparrow   &\geq \max \pbrk*{ v_{j-1/2}^+ , v_{j+1/2}^- }.
	\label{eqn:field_identities_vpm_bounds}
\end{align}
\end{subequations}

We design a well-balanced self-monotone reconstruction of $h$ by taking two separate reconstructions, and perform a convex combination \cref{eqn:convex_combination}. The two reconstructions are in depth $h$ and surface elevation $\eta$, yielding gradients in cell $C_j$
\begin{subequations} \label{eqn:results_gradient_h_eta}	\begin{align}
	[h_x]_j^{h} &\eqdef \sigma_j^h
	&\text{where}&&
	\sigma_j^h &\eqdef \sigma(h_{j-1},h_{j},h_{j+1};\{x\}),
	\label{eqn:results_gradient_h}
	\\
	[h_x]_j^{\eta} &\eqdef \sigma_j^\eta - [b_x]_j
	&\text{where}&&
	\sigma_j^\eta &\eqdef \sigma(\eta_{j-1},\eta_{j},\eta_{j+1};\{x\})
	\label{eqn:results_gradient_eta}
\end{align}\end{subequations}
respectively, using that the bed elevation is a known function with values $b_{j+1/2}^\pm$, \ie we include the case of discontinuities at interfaces. The interface values are found from \cref{eqn:results_gradient_h_eta} using \cref{eqn:FV_recon} to produce $h_{j+1/2}^{h\pm}$ and $h_{j+1/2}^{\eta\pm}$. Employing the convex combination \cref{eqn:convex_combination} we obtain a piecewise linear function with gradient
\begin{equation} \label{eqn:reaults_gradient_convex}
	[h_x]_j \eqdef (1-\gamma_j) [h_x]_j^{h} + \gamma_j [h_x]_j^{\eta},
\end{equation}
and values at the cell interfaces
\begin{subequations}\begin{align}
	h_{j-1/2}^{+} &\eqdef h_j - \frac{\Delta x_j}{2} [h_x]_j = (1-\gamma_j) h_{j-1/2}^{h+} + \gamma_j h_{j-1/2}^{\eta+},
	\\
	h_{j+1/2}^{-} &\eqdef h_j + \frac{\Delta x_j}{2} [h_x]_j = (1-\gamma_j) h_{j+1/2}^{h-} + \gamma_j h_{j+1/2}^{\eta-}.
\end{align}\end{subequations}
We use $\gamma_j$ to transition from \cref{eqn:results_gradient_h} to \cref{eqn:results_gradient_eta} by shifting from $\gamma_j=0$ to $\gamma_j=1$. The values $h_{j+1/2}^\pm$ are our reconstructed depths.

\begin{theorem}\label{thm:results_h}	\begin{subequations}
	Suppose that
	\begin{align}
		\gamma_j(\xi_j) &=
		\begin{cases}
			0 & \textrm{if } \xi \leq 1,	\\
			G_j(\xi_j - 1)	& \textrm{if } 1 \leq \xi_j \leq \xi_j^{C},	\\
			1 & \textrm{if } \xi_j^{C} \leq \xi_j,
		\end{cases}
		&&\textrm{where}&
		\xi_j^{C}(\{x\}) &= 1 + \frac{1}{G_j(\{x\})},
	\end{align}
	and $0 < G_j \leq 1-\alpha_j^\uparrow$, $\xi_j = {h_j^\downarrow}/{\Delta b_j^\uparrow}$ ($=+\infty$ for $\Delta b_j^\uparrow=0$) where
	\begin{align}\label{eqn:results_dbmax}
		\begin{split}
			\Delta b_j^\uparrow &\geq \max \left[ \abs{\Delta b_j/2-\alpha_{j-1/2}^{+}\Delta b_{j-1/2}},\abs{\Delta b_j/2}, 
		\right. \\ & \qquad \left.
			\abs{\Delta b_j/2 - \alpha_{j}(b_{j+1}-b_{j-1})},\abs{\Delta b_j/2 - \alpha_{j+1/2}^{-}\Delta b_{j+1/2}} \right],
		\end{split}
	\end{align}
	and $\Delta b_j^\uparrow$ is independent of $\{h\}$, then the reconstruction is self-monotone, the rate of decrease with the neighbouring cell values is bounded by
	\begin{align}\label{eqn:results_h_neigh}
		\pdv{h_{j-1/2}^+}{h_{j-1}} &\geq -G_j \alpha_{j-1/2}^{+} \geq -\frac{1}{4},
		&%\textrm{and}&&
		\pdv{h_{j+1/2}^-}{h_{j+1}} &\geq -G_j \alpha_{j+1/2}^{-} \geq -\frac{1}{4},
	\end{align}
	and the depth itself is bounded as
	\begin{equation}	\label{eqn:results_h_bound}	\begin{aligned}
		&&
		\ppar*{ 1-\frac{1}{\xi_j^{C}} } h_j^\downarrow &\leq h_{j-1/2}^+ \leq h_j^\uparrow + \frac{h_j^\downarrow}{\xi_j^{C}} ,
		\\&\text{and}&
		\ppar*{ 1-\frac{1}{\xi_j^{C}} } h_j^\downarrow &\leq h_{j+1/2}^- \leq h_j^\uparrow + \frac{h_j^\downarrow}{\xi_j^{C}}.
	\end{aligned}\end{equation}
\end{subequations}\end{theorem}
The absence of neighbour-monotonicity is discussed in \cref{sec:comparison,sec:supp_reconH}. For implementation, we note that the expression for $\gamma_j$ is continuous and piecewise linear in $\xi_j$, and is comparable in computational complexity to the $\minmod$ function. The value for $G_j$ chosen may be any value in the specified range, but typically taking $G_j = 1-\alpha_j^\uparrow$ should be the best choice. 
The expression for $\Delta b_j^\uparrow$ is given as an inequality to allow for reconstruction in $h$ is other situations, \eg by \cref{eqn:thin_film} we should reconstruct in $h$ when the local Froude number $u/\sqrt{gh}$ is large. We take $\Delta b_j^\uparrow$ as the larger of \cref{eqn:results_dbmax} and $B_j$, where $B_j$ is some function independent of $\{h\}$ \ie
\begin{equation}\label{eqn:results_db}\begin{split}
	\Delta b_j^\uparrow &= \max \left[ \abs*{\Delta b_j/2-\alpha_{j-1/2}^{+}\Delta b_{j-1/2}},\abs*{\Delta b_j/2}, 
	\right. \\ & \quad \left.
	\abs*{\Delta b_j/2 - \alpha_{j}(b_{j+1}-b_{j-1})},\abs*{\Delta b_j/2 - \alpha_{j+1/2}^{-}\Delta b_{j+1/2}} , B_j \right],
\end{split}\end{equation}
when $h_j^\downarrow \leq B_j$ we reconstruct in $h$. To capture `thin film' states take
\begin{equation}	\label{eqn:B_fastfluid}
	B_j = B(q_j) = \ppar*{\frac{q_j^2}{\Fro^2 g}}^{1/3}
\end{equation}
where $\Fro > 0$ is the reference Froude number. Thus when $\abs*{q_j}/[g (h_j^\downarrow)^3]^{1/2} \geq \Fro$ then the reconstruction is in $h$, and only when $\abs*{q_j}/[g (h_j^\downarrow)^3]^{1/2} \leq \Fro/(\xi_j^{C})^{3/2}$ is it possible to have a reconstruction in $\eta$; typically $\Fro \approx (\xi_j^{C})^{3/2}$ should be appropriate.

We next consider the velocity, $u$, which is computed using $u_j \eqdef q_j / h_j$ and $u_{j+1/2}^\pm \eqdef q_{j+1/2}^\pm/h_{j+1/2}^\pm$. We suppress the reconstruction using \cref{eqn:flux_suppress}, which yields
\begin{equation} \label{eqn:results_fluxrecon}
	[q_x]_j = \kappa_j \sigma (q_{j-1} , q_j , q_{j+1} ; \{x\}),
	\;\;\;\text{where}\;\;\;
	\kappa_j \eqdef \min \pbrk*{ 1 ,  \frac{K_{j-1/2}^+ h_j}{h_{j-1}} , \frac{K_{j+1/2}^- h_j}{h_{j+1}} },
\end{equation}
(to evaluate the ratios in the expression for $\kappa_j$ use $1/0 = +\infty$, no restriction on reconstruction, $0/0=0$, zero gradient in dry cell) where $K_{j+1/2}^{\pm} : \{x\} \rightarrow (0,\infty)$ determine how much deeper than the current cell a neighbour has to be for the gradient to be reduced. We recommend $K_{j+1/2}^{\pm} > 1  + \order{\Delta x_j}$ so that the gradient is only suppressed when $h_j$ is greater than $h_{j \pm 1}$, and in an implementation a value of $K_{j+1/2}^{\pm} \approx 100$ has been used successfully.
\begin{theorem}	\label{thm:results_u}
	With depth reconstructed as in \cref{thm:results_h} and $h_j>0$, and the volume flux reconstructed using \cref{eqn:results_fluxrecon}, the velocities have bounds
	\begin{equation}\begin{aligned}
		&&
		\abs*{u_{j-1/2}^+} &\leq \frac{G_j + 1}{1 - \alpha_j^{\uparrow}} \ppar*{ \abs*{u_j} + K_{j-1/2}^+ \alpha_{j-1/2}^{+} \abs*{u_{j-1}} },
		\\&\text{and}&
		\abs*{u_{j+1/2}^-} &\leq \frac{G_j + 1}{1 - \alpha_j^{\uparrow}} \ppar*{ \abs*{u_j} + K_{j+1/2}^- \alpha_{j+1/2}^{-} \abs*{u_{j+1}} }.
	\end{aligned}\end{equation}
\end{theorem}

Our approach is applicable to more sophisticated systems of equations, of which we discuss two. Firstly, for flow along a channel of varying width $w(x)$ a known function, the governing system is
\begin{subequations}	\label{eqn:SW_width_Full}	\begin{align}
	\pdv{wh}{t} + \pdv{}{x}(uwh) &= 0,		\label{eqn:SW_width_Full_h}	\\
	\pdv{}{t}(uwh) + \pdv{}{x} \ppar*{u^2 w h + \frac{g w h^2}{2}} &= - g w h\dv{b}{x} + \frac{g h^2}{2} \dv{w}{x}.			\label{eqn:SW_width_Full_u}
\end{align}\end{subequations}
This has steady states of constant $\tilde{q} \eqdef uwh$ and $E$, thus for well-balancing the reconstruction should still be in $\eta$ for deep areas and $h$ for shallow. This system is of interest because we must obtain a reconstruction for $A \eqdef wh$, a product of two fields for which gradients are first deduced independently, which is a non-trivial process by the discussion in \cref{sec:reconWH}. We derive the following result.
\begin{theorem}	\label{thm:results_wh}	\begin{subequations}
	Suppose that $[h_x]_j$ is computed using \cref{eqn:reaults_gradient_convex} with $\gamma_j$ from \cref{thm:results_h}, and that
	\begin{equation}\label{eqn:results_wh_grad}
		[A_x]_j 
		= [h_x]_j \ppar*{ w_j - \frac{\Delta x_j}{2} \abs*{[w_x]_j} } +  h_j [w_x]_j
		= [h_x]_j w_j^\downarrow +  h_j [w_x]_j
	\end{equation}
	with the reconstructed $A_{j+1/2}^\pm$ computed using \cref{eqn:FV_recon} and the reconstructed depth found by $h_{j+1/2}^\pm = A_{j+1/2}^\pm/w_{j+1/2}^\pm$, then the reconstructed depths are self-monotone and satisfy bounds \cref{eqn:results_h_neigh} and \cref{eqn:results_h_bound}, and if $\Delta w_j = \order{\Delta x_j}$ then the lake at rest state is accurate to $\sup \ppar*{ h_{j+1/2}^{-} + b_{j+1/2}^{-} - \bar{\eta} } = \order{(\Delta x_j)^2} = \ppar*{ h_{j-1/2}^{+} + b_{j-1/2}^{+} - \bar{\eta} }$.
	
	Moreover, if the flux is reconstructed as
	\begin{equation}
		[\tilde{q}_x]_j = \kappa_j \sigma (\tilde{q}_{j-1} , \tilde{q}_j , \tilde{q}_{j+1} ; \{x\}),
		\:\:\text{where}\:\:
		\kappa_j \eqdef \min \!\pbrk*{ 1 ,  \frac{K_{j-1/2}^+ A_j}{A_{j-1}} , \frac{K_{j+1/2}^- A_j}{A_{j+1}} }\!,
	\end{equation}
	where $K_{j+1/2}^\pm : \{x\} \rightarrow (0,\infty)$, then
	\begin{equation}\begin{aligned}
		&&
		\abs*{u_{j-1/2}^+} &\leq \frac{w_j}{w_j^\downarrow}\frac{G_j + 1}{1 - \alpha_j^{\uparrow}} \ppar*{ \abs*{u_j} + K_j \alpha_{j-1/2}^{+} \abs*{u_{j-1}} },
		\\ \text{and}&&
		\abs*{u_{j+1/2}^-} &\leq \frac{w_j}{w_j^\downarrow}\frac{G_j + 1}{1 - \alpha_j^{\uparrow}} \ppar*{ \abs*{u_j} + K_j \alpha_{j+1/2}^{-} \abs*{u_{j+1}} }.
	\end{aligned}\end{equation}
\end{subequations}\end{theorem}
The depth reconstruction is positive, self-monotone, and well-balanced, and we have bounded velocities everywhere where the width is not going to zero. If for some cell $w_j^\downarrow = 0$ then the velocities are unbounded. However, if say $w_{j-1/2}^+ = 0$ then we know that $\tilde{q}_{j-1/2}^+ = 0$, thus for this cell we reconstruct as $[\tilde{q}_x]_{j} = \tilde{q}_j 2 / \Delta x_j$, and in this cell $h$ is constant and $w$ and $q$ are linear, therefore $u$ is constant.

We discuss now a current driven by a density difference that varies depending on a particle concentration $\phi(x,t)$, which satisfies (\eg \cite{bk_Ungarish_GCI})
\begin{subequations}\label{eqn:EXT_SW_Full}
\begin{align}
	\pdv{h}{t} + \pdv{}{x}(uh) &= 0,		\label{eqn:EXT_SW_Full_h}	\\
	\pdv{}{t}(\phi h) + \pdv{}{x}(u \phi h) &= -v_s \phi,			\label{eqn:EXT_SW_Full_p}	\\
	\pdv{}{t}(uh) + \pdv{}{x} \ppar*{ u^2 h + \frac{g'_T h^2}{2} } &= - g'_T h\pdv{b}{x},			\label{eqn:EXT_SW_Full_u}
\end{align}
\end{subequations}
where $v_s \geq 0$ is the settling velocity, and $g'_T(x,t) \eqdef g'_p\phi(x,t) + g'_a$ is the total reduced gravity with $g'_p > 0$ and $g'_a \geq 0$. When $v_s=0$ this equation has steady states where $\phi$, $q$ and $E$ are constant, where $E$ is as in \cref{eqn:SW_steady_moving} with $g$ replaced by $g'_T$. The system exhibits `lake at rest' and `thin film' steady states, and to resolve these we use
\begin{equation}	\label{eqn:B_fastfluid_concen}
	B_j = B(\phi_j , q_j) = \ppar*{\frac{q_j^2}{\Fro^2 (g'_p \phi_j + g'_a)}}^{\mathrlap{1/3}}.
\end{equation}
This case is complicated for similar reasons to \cref{eqn:SW_width_Full}, and in addition the reconstruction in $h$ is dependent on $\{\phi\}$ by $B_j$. Despite these complications we achieve the following result.
\begin{theorem} \label{thm:results_phi} \begin{subequations}
	Suppose that $h$ is reconstructed as in \cref{thm:results_h} using \cref{eqn:results_db,eqn:B_fastfluid_concen}, and that
	\begin{equation}\label{eqn:results_phi_grad}
		[\Phi_x]_j = [\phi_x]_j \ppar*{ h_j - \frac{\Delta x_j}{2} \abs{[h_x]_j} } + \phi_j [h_x]_j
	\end{equation}
	where $\Phi \eqdef \phi h$, $[\phi_x]_j = \sigma(\phi_{j-1},\phi_{j},\phi_{j+1}; \{x\})$, $\phi_j = \Phi_j / h_j$, $\Phi_{j+1/2}^\pm$ are reconstructed as in \cref{eqn:FV_recon} , and $\phi_{j+1/2}^\pm = \Phi_{j+1/2}^+ / h_{j+1/2}^+$, then the values $\phi_{j+1/2}^\pm$ satisfy the self and neighbour-monotonicity properties, and are bounded as
	\begin{align}
		\phi_j^\downarrow &\leq \phi_{j-1/2}^+ \leq \phi_j^\uparrow,
		&
		\phi_j^\downarrow &\leq \phi_{j+1/2}^- \leq \phi_j^\uparrow.
	\end{align}
\end{subequations}\end{theorem}
The gradients \cref{eqn:results_wh_grad} and \cref{eqn:results_phi_grad} have the same structure following the discussion in \cref{sec:reconWH}. When $v_s>0$ we expect that the reconstruction \cref{eqn:results_phi_grad} is still appropriate, though the temporal evolution of the system is no longer straightforward because both the sink $-v_s \phi$ and the flux $u \phi h$ may act to remove particles from the cells, in which case $\phi$ can become negative after a time-step despite the reconstruction being positive. Positivity preserving Euler time-steps are presented in \cref{sec:postimestep}, generalizing of the method from \cite{ar_Bollermann_2011,ar_Bollermann_2013} and being compatible with the Runge-Kutta schemes in \cite{ar_Shu_1988}.
\section{Comparison} \label{sec:comparison}

We compare our depth reconstruction from \cref{thm:results_h} with others in the literature for the transition between `lake at rest' and `thin film' states plotted in \cref{fig:PWB_draining_flow}. We use a coarse grid with $x_{j+1/2} = j+ \frac*{1}{2}$ which is equivalent to a zoomed in view of a much finer grid. For smooth bed functions the gradient will be constant local to the transition, so we take $b_{j+1/2}^{\pm} = x_{j+1/2}$. The depth field used is $h_j = 1-b_j$ for $j \leq -2$ so that the surface is constant at $\eta = 1$, $h_{-1}=3$, $h_j = 1$ for $j \geq 0$ so that the depth is constant, and $h_0$ will take a range of values in each plot $0 \leq h_0 \leq 2$. This means that both the $h$ and $\eta$ fields are at a maxima in cell $C_{-1}$ so that $[h_x]_{-1}^h = 0$ and $[h_x]_{-1}^\eta = 1$, and $0 \leq \xi_{-1} \leq 6$, which allows us to investigate the monotonicity properties and lower bounds for each reconstruction. To neglect Froude number considerations we take \cref{eqn:results_dbmax} as equality.

\begin{figure}[tp!]
	\begin{center}
		\includegraphics{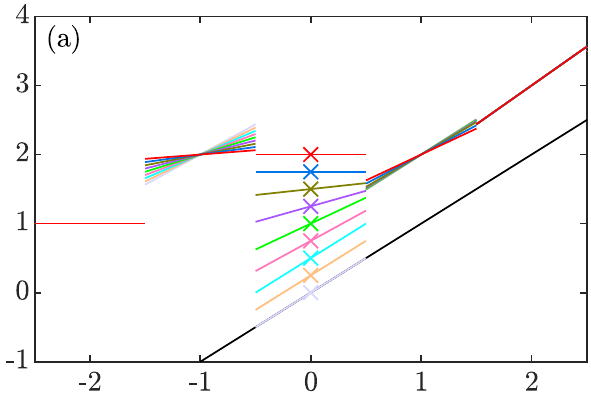}
		\includegraphics{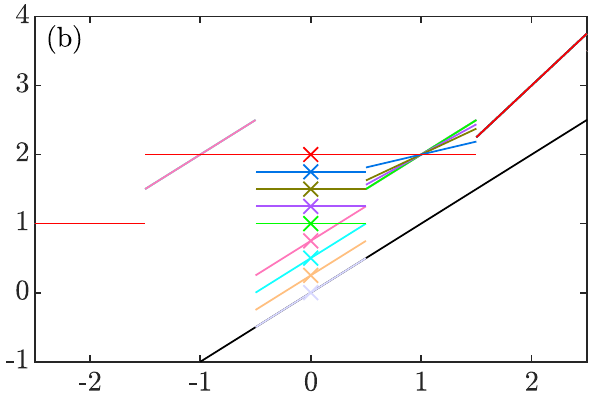}
		\includegraphics{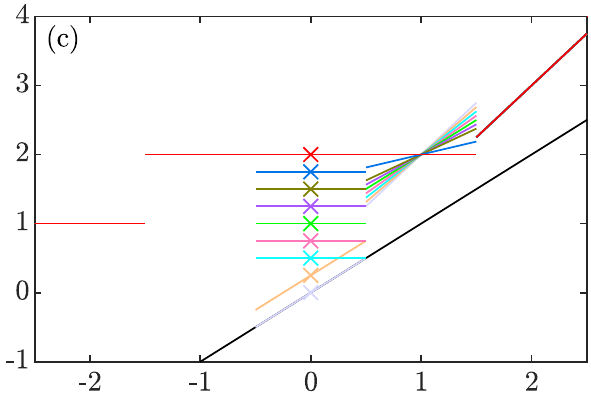}
		\includegraphics{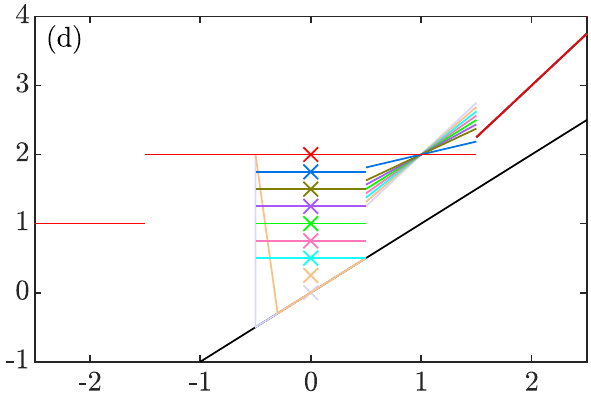}%
	\end{center}
	\caption{Reconstructed depths $h$ as a function of $x$, (a) uses our scheme, and (b), (c), and (d) use those of \cite{ar_Kurganov_2002}, \cite{ar_Chertock_2015}, and \cite{ar_Bollermann_2013} respectively. For each we consider multiple values of $h_0$, specifically $h_0=0,0.25,\ldots,2$, where each value is shown with a coloured cross at $x=0$. The reconstruction over all $x$ for this value of $h_0$ is plotted in this same colour. The gradient in $C_{-2}$ and $C_{2}$ is independent of $h_0$ for all schemes used, thus the plotted lines overlie and only the red can be seen. For (b), in $C_{-1}$ the reconstructed depths overlie for all $h_0\leq0.75$ (pink) and $h_0 > 0.75$ (red), and in $C_{1}$ they overlie for $h_0 \leq 1$ (green). For (c) and (d), in $C_{-1}$ all reconstructed depths overlie (red).}
	\label{fig:PWB_slopelim}
\end{figure}

Plots of the reconstructions are presented in \cref{fig:PWB_slopelim} for $\alpha_{j+1/2}^{\pm} = 3/4$, $\alpha_j = 1/4$, $G_j = 1/4$. In (a) the reconstruction is performed using the scheme presented here, in (b) by that of Kurganov and Levy \cite{ar_Kurganov_2002} with threshold $K=3/4$, in (c) by that of Chertock \etal \cite{ar_Chertock_2015}, and in (d) by that of Bollermann \etal \cite{ar_Bollermann_2013}. In (a) and (b) the reconstructions are finite for finite fluid depth, whereas (c) has $h_{1/2}^-=0$ when $h_0=1/2$, and in (d) $h_{1/2}^-=0$ for all of $0 \leq h_0 \leq 1/2$. Our reconstruction (a) is the only one which is self-monotone, or even has reconstructed values that are continuous in $h_0$. (b) has a discontinuous decrease in $h_{1/2}^-$ as $h_0$ increases past $K=3/4$, (c) has a similar decrease as $h_0$ increases past $b_{1/2}=1/2$, while (d) has a discontinuous decrease in $h_{-1/2}^+$ as $h_0$ increases past $b_{1/2}=1/2$. However, (c) and (d) are neighbour-monotone (this can be shown analytically for Chertock's in all cases, for Bollermann's it is possible to find edge cases where it is not), while for (b) there is a discontinuous decrease in $h_{-1/2}^-$ as $h_0$ increases past $K=3/4$, and for ours (a) $h_{-1/2}^-$ decreases as $h_0$ increases from $0$ to $7/3$ (at which point $\xi_{-1} = \xi_{-1}^C$) with ${\partial h_{-1/2}^-}/{\partial h_0} = -3/16$, the lower bound in \cref{thm:results_h}.

These observations give us confidence in our reconstruction.
\section{Monotone reconstruction for scalar problems} \label{sec:monotone}

While monotonicity (\cref{def:monotone}) is a reasonable requirement regardless of the type of problem considered, it is related to some other properties which we now discuss for the scalar function $v$ satisfying a scalar conservation law \cref{eqn:FV_conservation} with flux $F = f(v)$ and source $\Psi = 0$.

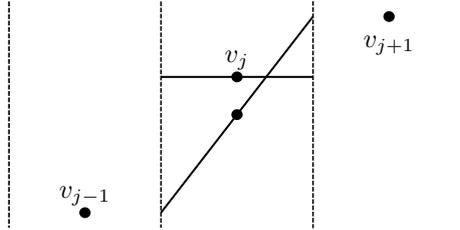
\begin{figure}[tp!]
	\centering
%% Accompanies image file 'NonSelfMonotone.pdf' (pdf, eps, ps)
\begingroup%
  \makeatletter%
  \providecommand\color[2][]{%
    \errmessage{(Inkscape) Color is used for the text in Inkscape, but the package 'color.sty' is not loaded}%
    \renewcommand\color[2][]{}%
  }%
  \providecommand\transparent[1]{%
    \errmessage{(Inkscape) Transparency is used (non-zero) for the text in Inkscape, but the package 'transparent.sty' is not loaded}%
    \renewcommand\transparent[1]{}%
  }%
  \providecommand\rotatebox[2]{#2}%
  \newcommand*\fsize{\dimexpr\f@size pt\relax}%
  \newcommand*\lineheight[1]{\fontsize{\fsize}{#1\fsize}\selectfont}%
  \ifx\svgwidth\undefined%
    \setlength{\unitlength}{198.42519685bp}%
    \ifx\svgscale\undefined%
      \relax%
    \else%
      \setlength{\unitlength}{\unitlength * \real{\svgscale}}%
    \fi%
  \else%
    \setlength{\unitlength}{\svgwidth}%
  \fi%
  \global\let\svgwidth\undefined%
  \global\let\svgscale\undefined%
  \makeatother%
  \begin{picture}(1,0.42857143)%
    \lineheight{1}%
    \setlength\tabcolsep{0pt}%
    \put(0,0){\includegraphics[width=\unitlength,page=1]{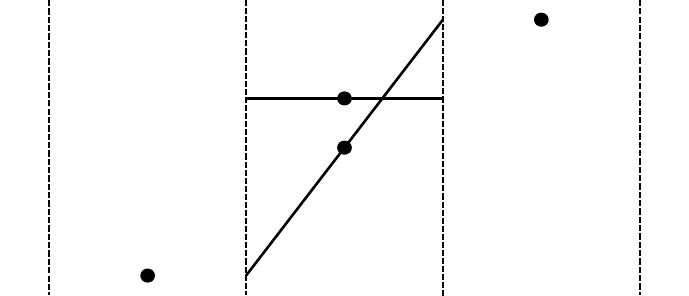}}%
    \put(0.21428571,0.05714286){\color[rgb]{0,0,0}\makebox(0,0)[t]{\lineheight{1.25}\smash{\begin{tabular}[t]{c}$v_{j-1}$\end{tabular}}}}%
    \put(0.78571429,0.34285714){\color[rgb]{0,0,0}\makebox(0,0)[t]{\lineheight{1.25}\smash{\begin{tabular}[t]{c}$v_{j+1}$\end{tabular}}}}%
    \put(0.5,0.31428571){\color[rgb]{0,0,0}\makebox(0,0)[t]{\lineheight{1.25}\smash{\begin{tabular}[t]{c}$v_{j}$\end{tabular}}}}%
  \end{picture}%
\endgroup%

	\caption{An example reconstruction $\hat{v}$ for cell $C_j$ where we consider two possible values of $v_j$. Increasing $v_j$ decreases $v_{j+1/2}^-$, thus it is not self-monotone.}
	\label{fig:PWB_non_self_monotone}
\end{figure}

We begin by discussing the use of the reconstruction in \cref{fig:PWB_non_self_monotone} in an upwind scheme for the advection equation. Increasing $v_j$ causes $v_{j+1/2}^{-}$ to decrease, decreasing the outflow causing $v_j$ to increase further as time passes until $v_{j+1/2}^{-}$ exceeds its initial value. When trying to decrease $v_j$ this feedback causes $v_j$ to decrease until $v_{j+1/2}^{-}$ is lower than its initial value. We propose the new requirements
\begin{align}	\label{eqn:self_neighbour_flux}
	\pdv{\tilde{f}}{v^-} (v^-,v^+) \cdot \pdv{f}{v} (v^-) &\geq 0,
	&\text{and}&&
	\pdv{\tilde{f}}{v^+} (v^-,v^+) \cdot \pdv{f}{v} (v^+) &\geq 0
\end{align}
for all $v^- , v^+$ on an interval where $f$ is monotone, where $f_{j+1/2} = \tilde{f}(v_{j+1/2}^{-},v_{j+1/2}^{+})$. In \cref{eqn:self_neighbour_flux} and going forward we make requirements for the full partial derivative, but equivalent requirements for the left and right partial derivatives can be made permitting discontinuous gradients at discrete points. Next we observe that, by self-monotonicity
\begin{subequations} \label{eqn:self_neighbour_gradient}
\begin{align}\label{eqn:self_gradient}
	\pdv{v_{j-1/2}^{+}}{v_{j}}  &\geq 0,
	&\textrm{and} &&
	\pdv{v_{j+1/2}^{-}}{v_{j}}  &\geq 0, 
\end{align}
and by neighbour-monotonicity
\begin{align}\label{eqn:neighbour_gradient}
	\pdv{v_{j-1/2}^{+}}{v_{j-1}}  &\geq 0,
	&\textrm{and} &&
	\pdv{v_{j+1/2}^{-}}{v_{j+1}}  &\geq 0.
\end{align}
\end{subequations}
From the requirements \cref{eqn:self_neighbour_flux} and \cref{eqn:self_neighbour_gradient}, in regions where the sign of $\pdv*{f}{v}$ is the same for $v_{j+1/2}^{+}$ and $v_{j+1/2}^{-}$, we have
\begin{align} \label{eqn:CDP_scheme}
	\pdv{f_{j+1/2}}{v_{j}} \cdot \pdv{f}{v} (v_{j+1/2}^{\pm}) &\geq 0,
	&\text{and}&&
	\pdv{f_{j+1/2}}{v_{j+1}} \cdot \pdv{f}{v} (v_{j+1/2}^{\pm}) &\geq 0.
\end{align}
This condition prevents the unstable behaviour discussed. We can interpret \cref{eqn:CDP_scheme} as a condition on motion of characteristics, if both $v_j$ and $v_{j+1}$ lead to characteristic motion in one direction then $f_{j+1/2}$ should also.
\begin{definition}
	A numerical flux is termed \emph{CDP (Characteristic Direction Preserving)} if it satisfies \cref{eqn:self_neighbour_flux}, a reconstruction is \emph{CDP} if it satisfies \cref{eqn:self_gradient} (a weaker condition than self-monotone) and \cref{eqn:neighbour_gradient}, and a scheme is termed \emph{CDP} if it satisfies \cref{eqn:CDP_scheme} when $v_{j+1/2}^{-}$ and $v_{j+1/2}^{+}$ are on an interval over which $f$ is monotone.
\end{definition}
\begin{lemma}	\label{thm:MUSCL_CDP}
	If the numerical flux and the reconstruction are CDP then the numerical scheme is CDP.
\end{lemma}
%\begin{proof}
%	Suppose the numerical flux and the reconstruction are CDP. On an interval where $f$ is monotone we have for $k \in \qty{j , j+1}$
%	\begin{equation*}
%		\pdv{f_{j+1/2}}{v_{k}} \cdot \pdv{f}{v} \qty(v_{j+1/2}^{\pm})
%		=
%		\qty[
%			\pdv{f_{j+1/2}}{v_{j+1/2}^{+}} \cdot \pdv{v_{j+1/2}^{+}}{v_k} + 
%			\pdv{f_{j+1/2}}{v_{j+1/2}^{-}} \cdot \pdv{v_{j+1/2}^{-}}{v_k}
%		]
%		\cdot \pdv{f}{v} \qty(v_{j+1/2}^{\pm})
%		\geq 0.
%	\end{equation*}
%	Thus the scheme is CDP.
%\end{proof}

While the definition of CDP schemes is new, the schemes themselves are not. Indeed any flux derived from consideration of the structure of the local Riemann problem should satisfy \cref{eqn:self_neighbour_flux}. For example exact Godunov (\eg \cite{bk_Leveque_FVM}), Local Lax-Friedrichs (also called Rusanov \cite{ar_Rusanov_1962} or Central \cite{ar_Kurganov_2000}), and Central-Upwind \cite{ar_Kurganov_2001} schemes. CDP reconstructions are not new either, for example the minmod \cite{ar_Tadmor_1988}, superbee \cite{ar_Roe_1985}, and MC \cite{ar_vanLeer_1976} slope limiters.

For inhomogeneous systems of equations such as \cref{eqn:INTRO_SW_Full} an appropriate extension of CDP is not clear, because one characteristic family may be travelling in one direction whilst another family is travelling in the other. For the purposes of this paper we enforce self-monotonicity for the well-balanced variables $W$. As per the discussion in the introduction we take $W = (h,q)^T$ (or equivalently $W=(\eta,q)^T$). Thus for self-monotonicity require
\begin{subequations} \label{eqn:self_neighbour_gradient_system}
	\begin{align}
	\pdv{W_{j-1/2}^{m+}}{W_{j}^{m}}  &\geq 0,
	&\textrm{and} &&
	\pdv{W_{j+1/2}^{m-}}{W_{j}^{m}}  &\geq 0, 
	\end{align}
and for neighbour-monotonicity 
	\begin{align}
	\pdv{W_{j-1/2}^{m+}}{W_{j-1}^{m}}  &\geq 0,
	&\textrm{and} &&
	\pdv{W_{j+1/2}^{m-}}{W_{j+1}^{m}}  &\geq 0,
	\end{align}
\end{subequations}
where $W^{m}$ is the $m^\textrm{th}$ component of $W$, and $W_{j+1/2}^{\pm}$ the transformed $Q_{j+1/2}^{\pm}$.
\section{Reconstruction of depth} \label{sec:reconH}

Here we prove \cref{thm:results_h}, beginning with self-monotonicity, for which we use that (after multiple applications of the chain rule)
\begin{subequations}\begin{align}
	\pdv{h_{j+1/2}^-}{h_j} &= \pdv{h_j^\downarrow}{h_j} \cdot \ppar*{ S_j - \pdv{\gamma_j}{\xi_j} R_j },	\label{eqn:reconH_selfmon_driv}
	\qquad \text{where}
	\\
	R_j &\eqdef \frac{h_{j+1/2}^{h-} - h_{j+1/2}^{\eta-}}{\Delta b_j^\uparrow}
	= \frac{\Delta b_j/2 + (\Delta x_j/2) (\sigma_j^h - \sigma_j^\eta)}{\Delta b_j^\uparrow},
	\\
	S_j &\eqdef \ppar*{\pdv{h_j^\downarrow}{h_j}}^{-1} \cdot \ppar*{ (1-\gamma_j) \pdv{h_{j+1/2}^{h-}}{h_j} + \gamma_j \pdv{h_{j+1/2}^{\eta-}}{h_j} }.
\end{align}\end{subequations}
Note that $\pdv*{h_j^\downarrow}{h_j} > 0$ because $\alpha_j^\uparrow < 1$. Proving results for $h_{j+1/2}^-$ is sufficient as it covers the case of $h_{j-1/2}^+$ by reflection in $x$. Next we bound $R_j$ and $S_j$.
\begin{lemma} \label{thm:reconH_bound_R}
	$-1 \leq R_{j} \leq 1$.
\end{lemma}
\begin{proof}\begin{subequations}
	We consider first when $h_{j-1} < h_{j} < h_{j+1}$ and $\eta_{j-1} < \eta_{j} < \eta_{j+1}$.
	\begin{itemize}[left=0mm]
		\item If $\sigma_j^{\eta}\Delta x_j/2 = \alpha_{j-1/2}^{+} \Delta \eta_{j-1/2}$ then by $\sigma_j^{h}\Delta x_j/2 \leq \alpha_{j-1/2}^{+} \Delta h_{j-1/2}$
		\begin{equation} \label{eqn:thm_reconH_bound_R_minmod_boundL}
		\eqnwide
			R_j
			\leq \frac{\Delta b_j/2 + \alpha_{j-1/2}^{+} (\Delta h_{j-1/2}-\Delta \eta_{j-1/2})}{\Delta b_j^\uparrow}
			= \frac{\Delta b_j/2 - \alpha_{j-1/2}^{+}\Delta b_{j-1/2}}{\Delta b_j^\uparrow}
			\leq 1.
		\end{equation}
		
		\item If $\sigma_j^{\eta}\Delta x_j/2 = \alpha_{j} (\eta_{j+1}-\eta_{j-1})$ then by $\sigma_j^{h}\Delta x_j/2 \leq \alpha_{j} (h_{j+1}-h_{j-1})$
		\begin{equation} \label{eqn:thm_reconH_bound_R_minmod_boundC} 
		\eqnwide
			R_j
			\leq \frac{\Delta b_j/2 + \alpha_{j} (h_{j+1}-h_{j-1}-\eta_{j+1}+\eta_{j-1})}{\Delta b_j^\uparrow}
			= \frac{\Delta b_j/2 - \alpha_{j} (b_{j+1} - b_{j-1})}{\Delta b_j^\uparrow}
			\leq 1.
		\end{equation}
		
		\item If $\sigma_j^{\eta}\Delta x_j/2 = \alpha_{j+1/2}^{-} \Delta \eta_{j+1/2}$ then by $\sigma_j^{h}\Delta x_j/2 \leq \alpha_{j+1/2}^{-} \Delta h_{j+1/2}$
		\begin{equation} \label{eqn:thm_reconH_bound_R_minmod_boundR}
		\eqnwide
			R_j
			\leq \frac{\Delta b_j/2 + \alpha_{j+1/2}^{-} (\Delta h_{j+1/2}-\Delta \eta_{j+1/2})}{\Delta b_j^\uparrow}
			= \frac{\Delta b_j/2 - \alpha_{j+1/2}^{-} \Delta b_{j+1/2}}{\Delta b_j^\uparrow}
			\leq 1.
		\end{equation}
	\end{itemize}
	Consider now $h_{j-1} > h_{j} > h_{j+1}$ and $\eta_{j-1} > \eta_{j} > \eta_{j+1}$.
	\begin{itemize}[left=0mm]
		\item If $\sigma_j^{h}\Delta x_j/2 = \alpha_{j-1/2}^{+} \Delta h_{j-1/2}$ then by $\sigma_j^{\eta}\Delta x_j/2 \geq \alpha_{j-1/2}^{+} \Delta \eta_{j-1/2}$ we obtain \cref{eqn:thm_reconH_bound_R_minmod_boundL}.
		
		\item If $\sigma_j^{h}\Delta x_j/2 =  \alpha_{j} (h_{j+1}-h_{j-1})$ then by $\sigma_j^{\eta}\Delta x_j/2 \geq \alpha_{j} (\eta_{j+1}-\eta_{j-1})$ we obtain \cref{eqn:thm_reconH_bound_R_minmod_boundC}.
		
		\item If $\sigma_j^{h}\Delta x_j/2 = \alpha_{j+1/2}^{-} \Delta h_{j+1/2}$ then by $\sigma_j^{\eta}\Delta x_j/2 \geq \alpha_{j+1/2}^{-} \Delta \eta_{j+1/2}$ we obtain \cref{eqn:thm_reconH_bound_R_minmod_boundR}.
	\end{itemize}
	Finally, the cases where the discretised variables are not monotone of the same sign.
	\begin{itemize}[left=0mm]
		\item If [$\Delta h_{j-1/2} \leq 0$ or $\Delta h_{j+1/2} \leq 0$] and [$\Delta \eta_{j-1/2} \geq 0$ or $\Delta \eta_{j+1/2} \geq 0$] then
		\begin{equation}
		\eqnwide
			R_j
			\leq \frac{\Delta b_j/2}{\Delta b_j^\uparrow}
			\leq 1.
		\end{equation}
		\item If $\Delta h_{j-1/2} \geq 0$ and $\Delta\eta_{j-1/2} \leq 0$ then we obtain \cref{eqn:thm_reconH_bound_R_minmod_boundL}.
		\item If $\Delta h_{j+1/2} \geq 0$ and $\Delta\eta_{j+1/2} \leq 0$ then we obtain \cref{eqn:thm_reconH_bound_R_minmod_boundR}.
	\end{itemize}
	The result $R_j \geq -1$ comes from reflection in $x$, under which $\Delta b$, $\sigma_j^h$ and $\sigma_j^\eta$ change sign and $\Delta x_j$ and $\Delta b_j^\uparrow$ do not.
\end{subequations}\end{proof}

\begin{lemma}	\label{thm:reconH_bound_S}
	$S_j \geq 1 - \alpha_{j+1/2}^-.$
\end{lemma}
\begin{proof}
	First observe that the minmod slope limiter \cref{eqn:minmod_recon_R} satisfies
	\begin{gather*}
		\pdv{\sigma}{v_j} \geq - \frac{2 \alpha_{j+1/2}^-}{\Delta x_j},
		\qquad\text{and}\qquad
		1 - \alpha_j^\uparrow \leq \pdv{h_j^\downarrow}{h_j} \leq 1,
	\shortintertext{thus}
		S_j \geq 1 \cdot \ppar*{ (1-\gamma_j)(1 - \alpha_{j+1/2}^-) + \gamma_j(1 - \alpha_{j+1/2}^-) } = 1 - \alpha_{j+1/2}^-.
	\end{gather*}
\end{proof}

\begin{lemma}	\label{eqn:reconH_selfmon}
	For $\gamma_j$ as in \cref{thm:results_h} the reconstruction is self monotone.
\end{lemma}
\begin{proof}
	Employing $0 \leq \pdv*{\gamma_j}{\xi_j} \leq G_j$, \cref{eqn:reconH_selfmon_driv,thm:reconH_bound_R,thm:reconH_bound_S} we obtain
	\begin{equation*}
		\pdv{h_{j+1/2}^-}{h_j} \geq \pdv{h_j^\downarrow}{h_j} \cdot \ppar*{ 1 - \alpha_{j+1/2}^- - G_j } \geq 0.
	\end{equation*}
\end{proof}

Next we prove the bounds on the depth and its derivatives, focussing on $h_{j+1/2}^-$.

\begin{lemma}	\label{thm:reconH_neigh_severe}
	The derivatives of the depth reconstruction have lower bounds
	\begin{align}
		\pdv{h_{j+1/2}^-}{h_{j+1}} &\geq -G_j \alpha_{j+1/2}^{-} \geq -\frac{1}{4}.
	\end{align}
\end{lemma}
\begin{proof}
	By symmetry we only need to prove for $h_{j+1/2}^-$, for which
	\begin{equation}
		\pdv{h_{j+1/2}^-}{h_{j+1}} = \ppar*{ (1-\gamma_j) \pdv{h_{j+1/2}^{h-}}{h_j} + \gamma_j \pdv{h_{j+1/2}^{\eta-}}{h_j} } - \pdv{h_j^\downarrow}{h_{j+1}} \pdv{\gamma_j}{\xi_j} R_j.
	\end{equation}
	This is positive unless $h_{j+1} \leq \min[h_{j},h_{j-1}]$ and $\pdv*{h_j^\downarrow}{h_{j+1}} = \alpha_{j+1/2}^-$, in which case
	\begin{equation}
		\pdv{h_{j+1/2}^-}{h_{j+1}} \geq -G_j \alpha_{j+1/2}^{-} \geq -(1-\alpha_{j}^{\uparrow}) \alpha_{j+1/2}^{-}.
	\end{equation}
	Using that $0 < \alpha_{j+1/2}^{-} \leq \alpha_{j}^{\uparrow} < 1$, this product takes its minimum when $\alpha_{j+1/2}^{-} = \alpha_{j}^{\uparrow} = 1/2$ with a value of $-1/4$.
\end{proof}

\begin{lemma}	\label{thm:reconH_h_bound}	\begin{subequations}
	The depth reconstruction has bounds
	\begin{equation*}
		\ppar*{ 1-\frac{1}{\xi_j^{C}} } h_j^\downarrow \leq h_{j+1/2}^- \leq h_j^\uparrow + \frac{h_j^\downarrow}{\xi_j^{C}} .
	\end{equation*}
\end{subequations}\end{lemma}
\begin{proof}
	Firstly we note that $h_j^\downarrow \leq h_{j+1/2}^{h-} \leq h_j^\uparrow$, and 
	\begin{subequations}\begin{align*}
		\begin{split}
			h_{j+1/2}^{\eta-}
			&\geq \min\pbrk*{h_j - \Delta b_j/2 + \alpha_{j+1/2}^{-}\Delta \eta_{j+1/2} , h_j - \Delta b_j/2}
			\\
			&\geq h_{j}^\downarrow - \Delta b_j^\uparrow
			= (1 - (\xi_j)^{-1})h_j^\downarrow,
		\end{split}
		\displaybreak[3]\\
		\begin{split}
			h_{j+1/2}^{\eta-}
			&\leq \max\pbrk*{h_j - \Delta b_j/2 + \alpha_{j+1/2}^{-}\Delta \eta_{j+1/2} , h_j - \Delta b_j/2}
			\\
			&\leq h_{j}^\uparrow + \Delta b_j^\uparrow
			= h_j^\uparrow + (\xi_j)^{-1} h_j^\downarrow.
		\end{split}
	\end{align*}\end{subequations}
	When $\xi_j \leq 1$, $h_{j+1/2}^- = h_{j+1/2}^{h-}$, and when $\xi_j \geq (\xi_j^{C})^{-1}$, $h_{j+1/2}^- = h_{j+1/2}^{\eta-}$, thus in these cases the bounds are already verified.
	When $1 \leq \xi \leq \xi_j^{C}$,
	\begin{subequations}\begin{align*}
		\begin{split}
			h_{j+1/2}^- &= (1-\gamma_j) h_{j+1/2}^{h-} + \gamma_j h_{j+1/2}^{\eta-}
			\geq (1-\gamma_j)h_j^\downarrow + \gamma_j (1 - (\xi_j)^{-1}) h_j^\downarrow
			\\&
			= (1 + G_j ((\xi_j)^{-1} - 1)) h_j^\downarrow
			\geq (1 + G_j ((\xi_j^C)^{-1} - 1)) h_j^\downarrow
			= (1 - (\xi_j^C)^{-1}) h_j^\downarrow,
		\end{split}
		\displaybreak[3]\\
		\begin{split}
			h_{j+1/2}^- &= (1-\gamma_j) h_{j+1/2}^{h-} + \gamma_j h_{j+1/2}^{\eta-}
			\leq (1-\gamma_j) h_j^\uparrow + \gamma_j (h_j^\uparrow + (\xi_j)^{-1} h_j^\downarrow)
			\\&
			= h_j^\uparrow + G_j (1-(\xi_j)^{-1}) h_j^\downarrow
			\leq h_j^\uparrow + G_j (1-(\xi_j^C)^{-1}) h_j^\downarrow
			= h_j^\uparrow + (\xi_j^C)^{-1} h_j^\downarrow.
		\end{split}
	\end{align*}\end{subequations}
\end{proof}
\section{Reconstruction of flux} \label{sec:recon_q}

In this section we prove \cref{thm:results_u}. First
\begin{align}\begin{split}
	\abs*{q_{j+1/2}^-} 
	&\leq \max \pbrk*{ \abs*{q_j} , \abs*{q_j + \kappa_j \alpha_{j+1/2}^{-} (q_{j+1}-q_j)}  }
	\\
	&\leq \max \pbrk*{ \abs*{u_j} h_j , (1 - \kappa_j \alpha_{j+1/2}^{-}) \abs*{u_j} h_j + \kappa_j \alpha_{j+1/2}^{-} \abs*{u_{j+1}} h_{j+1} }.
\end{split}\end{align}
Next we use that, no matter the discretized values of depth, $\kappa_j h_{j+1} \leq K_{j+1/2}^- h_j$, thus
\begin{align}\begin{split}
	\abs*{q_{j+1/2}^-} 
	&\leq h_j \max \pbrk*{ \abs*{u_j} , (1 - \kappa_j \alpha_{j+1/2}^{q-}) \abs*{u_j} + K_{j+1/2}^- \alpha_{j+1/2}^{q-} \abs*{u_{j+1}} }
	\\
	&\leq h_j \ppar*{ \abs*{u_j} + K_{j+1/2}^- \alpha_{j+1/2}^{q-} \abs*{u_{j+1}} }.
\end{split}\end{align}
To turn this into a bound on $\abs{u_{j+1/2}^-}$ we employ \cref{thm:reconH_h_bound}, specifically
\begin{align}
	\frac{h_j}{h_{j+1/2}^-} \leq (G_j+1) \frac{h_j}{h_j^\downarrow} \leq \frac{G_j + 1}{1 - \alpha_j^{\uparrow}}
\end{align}
where the second inequality is from consideration of the case $h_{j-1} = h_{j+1} = 0$, thus
\begin{align}
	\abs*{u_{j+1/2}^-} \leq \frac{G_j + 1}{1 - \alpha_j^{\uparrow}} \ppar*{ \abs*{u_j} + K_{j+1/2}^- \alpha_{j+1/2}^{-} \abs*{u_{j+1}} }
\end{align}
and employing symmetry under reflection in $x$ proves \cref{thm:results_u}
\section{Modified depth reconstruction for the inclusion of width} \label{sec:reconWH}

In this section we prove \cref{thm:results_wh}, and in addition give an overview of the problems that can occur when trying to reconstruct the gradient of a product of two functions that are initially reconstructed independently.

We being our discussion for the case of an arbitrary scalar field $v$, and once we have understood how to reconstruct $wv$ we discuss $wh$. The scalar field is calculated to have a gradient $[v_x]_j$, which yields reconstructed values $v_{j+1/2}^{v\pm}$ by \cref{eqn:FV_recon}; these are not the reconstructed values to be used in the scheme, that reconstruction will be in $wv$, but rather values used to aid discussion. Using $[v_x]_j$ we reconstruct $wv$, then compute the values $[wv]_{j+1/2}^\pm$ by \cref{eqn:FV_recon}, and finally $v_{j+1/2}^\pm \eqdef [wv]_{j+1/2}^\pm/w_{j+1/2}^\pm$, these are the values to be used. The reconstruction of $wv$ could be performed by simply multiplying the piecewise linear $w$ and $v$ to produce a quadratic, but this would then need to be modified so that the integral of the quadratic matches the cell averaged value $[wv]_j$, which encounters the same positivity issue we resolve here (see \cref{sec:suppreconWH_simp}). Instead, we use a piecewise linear reconstruction with gradient $[(wv)_x]_j$. We define three reconstructions indexed by the integer $r \in \{-1,0,1\}$ respectively representing left, centred, and right differences of the quadratic,
\begin{align} \label{eqn:reconWH_grad_stucture}
	[(wv)_x]_j^{(r)} &\eqdef [w_x]_j v_j + w_j [v_x]_j + r \frac{\Delta x_j}{2} [w_x]_j [v_x]_j.
\end{align}
The issue encountered with constructing this gradient is that, while $w$ and $v$ are individually constructed to be positive, there is no guarantee of a positive lower bound for $wv$. Indeed, taking $r=0$ gives an expression with the appearance of the product rule, but produces negative values when $\abs*{[w_x]_j/w_j + [v_x]_j/v_j} \geq 2/\Delta x_j$ (see \cref{sec:suppreconWH_simp}). 

We select a value of $r$ in each cell, $r_j$, which yields a reconstruction $[(w v)_x]_j \eqdef [w v_x]_j^{(r_j)}$. The reconstructed $v$ is given by
\begin{align}
	v_{j+1/2}^- &= v_j + \frac{w_j + r_j \frac{\Delta x_j}{2} [w_x]_j}{w_{j-1/2}^+} \frac{\Delta x_j}{2} [v_x]_j.
\end{align}
The difference between reconstructing in $wv$ and reconstructing in $v$ directly is that the effect of the gradient $[v_x]_j$ is multiplied by some ratio of widths $w$. Imposing that this ratio is at most one, then the reconstructed values $v_{j+1/2}^\pm$ immediately satisfy all bounds that may be deduced for $v_{j+1/2}^{v\pm}$. We take $r_j = - \sign(\Delta w_j)$, thus
\begin{gather} \label{eqn:reconWH_grad}
	[(wv)_x]_j = [w_x]_j v_j + w_j^\downarrow [v_x]_j, \\
	\begin{aligned}
		v_{j-1/2}^+ &= v_j - \frac{w_j^\downarrow}{w_{j-1/2}^+} \frac{\Delta x_j}{2} [v_x]_j,	&
		v_{j+1/2}^- &= v_j + \frac{w_j^\downarrow}{w_{j+1/2}^-} \frac{\Delta x_j}{2} [v_x]_j.
	\end{aligned}
\end{gather}
Bounding $v_{j+1/2}^-$ we obtain (see \cref{sec:suppreconWH_bound})
\begin{subequations} \label{eqn:reconWH_bounds}
\begin{align}
	\label{eqn:reconWH_value}
	\min \pbrk*{v_j , v_{j+1/2}^{v-}} &\leq v_{j+1/2}^{-} \leq \max \pbrk*{v_j , v_{j+1/2}^{v-}}
	\\
	\label{eqn:reconWH_selfmon}
	\min \pbrk*{ 1 , \pdv{v_{j+1/2}^{v-}}{v_j} } &\leq \pdv{v_{j+1/2}^{-}}{v_j} \leq \max \pbrk*{ 1 , \pdv{v_{j+1/2}^{v-}}{v_j} }
	\\
	\label{eqn:reconWH_neighmon}
	\min \pbrk*{ 0 , \pdv{v_{j+1/2}^{v-}}{v_{j+1}} } &\leq \pdv{v_{j+1/2}^{-}}{v_{j+1}} \leq \max \pbrk*{ 0 , \pdv{v_{j+1/2}^{v-}}{v_{j+1}} }
\end{align}
\end{subequations}
thus bounds on $v_{j+1/2}^{v-}$ and its derivatives extend straightforwardly to $v_{j+1/2}^{-}$.

Returning to the variable $h$ we have
\begin{gather}
	[A_x]_j = [w_x]_j h_j + w_j^\downarrow [h_x]_j, \\
	\begin{aligned}
		h_{j-1/2}^+ &= h_j - \frac{w_j^\downarrow}{w_{j-1/2}^+} \frac{\Delta x_j}{2} [h_x]_j,	&
		h_{j+1/2}^- &= h_j + \frac{w_j^\downarrow}{w_{j+1/2}^-} \frac{\Delta x_j}{2} [h_x]_j.
	\end{aligned}
\end{gather}
with $[h_x]_j$ from \cref{eqn:reaults_gradient_convex}. The bounds on $h_{j+1/2}^-$, $\pdv*{h_{j+1/2}^-}{h_j}$, and  $\pdv*{h_{j+1/2}^-}{h_{j+1}}$ from \cref{thm:results_h} still apply, ensuring positivity. In addition, for $\Delta w_j = \order{\Delta x_j}$ the reduction in the gradient is $\order{\Delta x_j}$, thus the modification in $h_{j+1/2}^-$ is $\order{\Delta x_j^2}$.

To deduce bounds on the velocity we first bound $A_{j+1/2}^-$ from below, where
\begin{equation}
	A_{j+1/2}^- = w_j \cdot \ppar*{h_j + \frac{\Delta x_j}{2}[h_x]_j} + \frac{\Delta x_j}{2} [w_x]_j \cdot
	\begin{cases}
		h_j - \frac{\Delta x_j}{2}[h_x]_j 
		&\textrm{for } [w_x]_j > 0,
		\\
		h_j + \frac{\Delta x_j}{2}[h_x]_j 
		&\textrm{for } [w_x]_j < 0.
	\end{cases}
\end{equation}
Using the bounds from \cref{thm:reconH_h_bound} we obtain
\begin{equation}
	A_{j+1/2}^- \geq \ppar*{1 - \frac{1}{\xi_j^C}} h_j^\downarrow w_j^\downarrow.
\end{equation}
The bounds on velocity then follow using the method from \cref{sec:recon_q}, which completes the proof of \cref{thm:results_wh}. 
\section{Reconstruction of concentration} \label{sec:recon_phi}

In this section we prove \cref{thm:results_phi}. Exactly as was the case in \cref{sec:reconWH}, from $[\phi_x]_j$ we construct a gradient $[\Phi_x]_j$, then the values $\Phi_{j+1/2}^\pm$ by \cref{eqn:FV_recon}, and finally $\phi_{j+1/2}^\pm \eqdef \Phi_{j+1/2}^\pm/h_{j+1/2}^\pm$, these are the values to be used. Thus we may use the results of \cref{sec:reconWH}, and immediately set
\begin{equation}
	[\phi h_x]_j \eqdef [\phi_x]_j \ppar*{ h_j - \frac{\Delta x_j}{2} \abs*{[h_x]_j} } + \phi_j [h_x]_j.
\end{equation}
as our reconstruction, and thus
\begin{align}
	\phi_j^\downarrow &\leq \phi_{j-1/2}^+ \leq \phi_j^\uparrow,
	&
	\phi_j^\downarrow &\leq \phi_{j+1/2}^- \leq \phi_j^\uparrow.
\end{align}
However, unlike the case of \cref{sec:reconWH}, the reconstruction of $h$ is dependent on $\{\phi\}$ though $B_j$ \cref{eqn:B_fastfluid_concen}. Because $B_j$ is independent of $\phi_{j \pm 1}$ we immediately get neighbour-monotonicity. Self-monotonicity is more challenging because $B_j$ is a function of $\phi_j$; we consider
\begin{align}
	\phi_{j+1/2}^- &= \phi_j + \frac{h_j - \frac{\Delta x_j}{2} \abs*{[h_x]_j}}{h_{j+1/2}^-} \frac{\Delta x_j}{2} [\phi_x]_j
\end{align}
$\phi_{j-1/2}^+$ similar by symmetry. For $[h_x]_j \leq 0$ we have $\phi_{j+1/2}^- = \phi_j + [\phi_x]_j {\Delta x_j}/{2}$, thus we have self-monotonicity. Also, when $\pdv*{[h_x]_j}{\phi_j} = 0$ then \cref{eqn:reconWH_selfmon} can be applied and, again, we have self-monotonicity. The only case left is $[h_x]_j \geq 0$, $1 \leq \xi_j \leq \xi_j^C$, and $\Delta b_j^\uparrow = B_j$, which we now consider. We first compute some results for the depth
\begin{gather*}
	\pdv{h_{j+1/2}^-}{\phi_j} 
	= \pdv{\gamma_j}{\xi_j} \cdot \pdv{\xi_j}{\Delta b_j^\uparrow} \cdot  \pdv{\Delta b_j^\uparrow}{\phi_j} \cdot \ppar*{ h_{j+1/2}^{\eta - } - h_{j+1/2}^{h - } }
	= G_j R_j \frac{h_{j}^\downarrow}{B_j} \pdv{B_j}{\phi_j},
\shortintertext{similarly}
	\pdv{h_{j-1/2}^+}{\phi_j} = -G_j R_j \frac{h_{j}^\downarrow}{B_j} \pdv{B_j}{\phi_j}.
\shortintertext{Thus}
\begin{split}
	\pdv{\phi_{j+1/2}^-}{\phi_j} &= 1 + \frac{\Delta x_j}{2} \pdv{}{\phi_j} \ppar*{\frac{h_{j-1/2}^+ [\phi_x]_j}{h_{j+1/2}^-}}	\\
	&= 1 + \frac{\Delta x_j}{2} \frac{h_{j-1/2}^+}{h_{j+1/2}^-} \pdv{[\phi_x]_j}{\phi_j} - \frac{\Delta x_j}{2} \frac{[\phi_x]_j}{h_{j+1/2}^-} \cdot G_j R_j \frac{h_{j}^\downarrow}{B_j} \pdv{B_j}{\phi_j} \cdot \ppar*{ 1 + \frac{h_{j-1/2}^+}{h_{j+1/2}^-} }.
\end{split}
\end{gather*}
Next we note that, for $B_j$ as given in \cref{eqn:B_fastfluid_concen},
\begin{equation}
	P_j \eqdef \sup{\abs*{\frac{\phi_j}{B_j} \pdv{B_j}{\phi_j} }} = \sup{\abs*{ -\frac{1}{3} \frac{g'_p \phi_j}{g'_p \phi_j + g'_a} }} = \frac{1}{3}.
\end{equation}
Using that $\pdv*{[\phi_x]_j}{\phi_j} \geq -2 \alpha_{j+1/2}^{-}/\Delta x_j$, $\abs*{[\phi_x]_j} \leq 2 \alpha_j^{\uparrow} \phi_j/\Delta x_j$, and $\abs*{R_j} \leq 1$ by \cref{thm:reconH_bound_R}, we obtain
\begin{align}
	\pdv{\phi_{j+1/2}^-}{\phi_j} &\geq 1 - \alpha_{j+1/2}^{-}  -  2 \alpha_j^{\uparrow} G_j P_j \geq \ppar*{1 - \alpha_{j+1/2}^{-}} \ppar*{ 1 - \frac{2}{3} \alpha_j^{\uparrow}} \geq 0
\end{align}
where we have used $[h_x]_j \geq 0$ to set all depth ratios to their maximal value of $1$. This proves \cref{thm:results_phi}.
\section{Conclusion}

We have established an approach (\cref{sec:overview}), valid for a wide range of systems, to design a well-balanced numerical scheme that is self-monotone (\cref{def:monotone}). This new property is in accordance with physical intuition, but also connects to stability considerations (\cref{sec:monotone}). The particular scheme that we develop for the shallow water equations is given in \cref{sec:results}, using the minmod reconstruction as a base, and also generalising to a range of other systems. The reconstruction that results is compared to other approaches in the literature (\cref{sec:comparison}), demonstrating that our reconstruction has many favourable properties. Comparisons for full numerical simulations are not presented here, and will appear in a future manuscript. The approach that we introduce (\cref{sec:overview}) is generalisable to a wide range of systems and may use an arbitrary reconstruction of any order as a base, meaning that high order well balanced schemes with good stability properties may be developed.
\section*{Acknowledgements}

The author would also like to thank A. J. Hogg for his constructive comments regarding drafts of this article.

\begin{myappendices}
	\section{Positivity preserving time-stepping with sources} \label{sec:postimestep}

%\begin{figure}
%	\centering
%	\input{Figures/PositivityPreserving.pdf_tex}
%	\caption{Schematic showing the ways that fluxes and sources can act on a cell, the fluxes are represented by horizontal arrows showing the flux to adjacent cells, the sources by vertical arrows. For the purposes of this schematic only we suppose that the sources extract material from the flow through the bed, as is the case for particle settling. The three diagrams illustrate: (a) the fluxes move material into the cell whilst the source drains it; (b) the fluxes move material across the cell whilst the source drains it; (c) the fluxes and source all act to drain the cell.}
%	\label{fig:PP_flux_sink}
%\end{figure}

After an Euler time-step of \cref{eqn:EXT_SW_Full} the values of $h$ and $\phi h$ in each cell must remain positive. So long as there are no sinks in the system the CFL condition with Courant number at most $1/2$ is sufficient to enforce this, as it gives $\abs*{u_{j+1/2}^-},\abs*{u_{j-1/2}^+} \leq \frac*{\Delta x_j}{2\Delta t}$. However, in \cref{eqn:EXT_SW_Full} the particle field has the possibility of sinks (and for models including other physical processes the depth field may have sinks too), and the CFL condition is insufficient to enforce positivity. To remedy this we propose a draining time technique similar to those in \cite[\S4]{ar_Bollermann_2011} and \cite[\S4]{ar_Bollermann_2013}, but modified to account for source terms. Unlike in a normal Euler time-step where the fluxes and sources act constantly for the full interval $\Delta t$, we instead suppose that they act constantly until all the material (fluid/particles) they have access to is drained and then cease for the remainder of the time step. The draining time for fluxes and sources is denoted by $\Delta t_j^{F,m}$ and $\Delta t_j^{\Psi,m}$ respectively for the $m^\textrm{th}$ field. This results in the fluxes and sources acting for proportions of the time-step
\begin{subequations}\begin{align}
	d_{j+1/2}^{F,m} &= \min \pbrk*{ \frac{\Delta t_{k}^{F,m}}{\Delta t} , 1 }
	\quad\textrm{where}\quad
	k = j+\frac{1-\sign(F_{j+1/2}^{m})}{2},
	\\
	d_j^{\Psi,m} &= \min \pbrk*{ \frac{\Delta t_j^{\Psi,m}}{\Delta t} , 1 }
\end{align}\end{subequations}
respectively, where the superscript $m$ denotes the $m^{\nth}$ component of the vector. The resultant fluxes and sources for use in the Euler time-step are
\begin{align}
	\tilde{F}_{j+1/2}^{m} &= d_{j+1/2}^{F,m} F_{j+1/2}^{m},
	&\textrm{and}&&
	\tilde{\Psi}_j^{m} &= d_j^{\Psi,m} \Psi_j^{m}.
\end{align}
When the source for the field is positive, $\Psi_j^{m} \geq 0$, and is adding more fluid/particles into the domain, then $\Delta t_j^{F,m} \geq \Delta t$ and $\Delta t_j^{\Psi,m} = \infty$ by CFL. When the source for the field is negative, $\Psi_j^{m} < 0$, then we use that the fluxes can only act until the material that was initially present in the cell is exhausted. Thus
\begin{equation}
	\Delta t_j^{F,m} = \frac{Q_j^{m} \Delta x_j}{\max\left(F_{j+1/2}^{m},0\right) + \max\left(-F_{j-1/2}^{m},0\right) + \max\left(-\Psi_j^{m}\Delta x_j,0\right)}.
\end{equation}
The sources not only have access to that material initially present in the cell, but also that material advected into the cell over the current time step. Therefore
\begin{equation}
	\Delta t_j^{\Psi,m} =
	\begin{cases}
		\infty
		& \textrm{for } \Psi_j^{m} \geq 0,
		\\
		\Delta t_j^{F,m} \! + \! \dfrac{\ppar*{\max\pbrk*{-\tilde{F}_{j+1/2}^{m},0} + \max\pbrk*{\tilde{F}_{j-1/2}^{m},0}} \Delta t}{-\Psi_j^{m} \Delta x_j}
		& \textrm{for } \Psi_j^{m} < 0.
	\end{cases}
\end{equation}

Following \cite{ar_Bollermann_2013}, fluxes in the momentum equation should also cease, only permitting the $u^2 h$ term to act whilst the cell is not drained of fluid. A similar modification ceases the flux of particles whilst there is no flux of fluid, thus
\begin{align}
	\tilde{d}_{j+1/2}^{F,2} &= \min\pbrk*{ d_{j+1/2}^{F,1}, d_{j+1/2}^{F,2} },
	&\textrm{and}&&
	\tilde{d}_j^{\Psi,2} &= \min\pbrk*{ d_j^{\Psi,1}, d_j^{\Psi,2} }.
\end{align}
The $g_T' h^2 /2$ flux term and $g_T' h \pdv*{b}{x}$ source term in the momentum equation should not be ceased as these balance each other.
\end{myappendices}

\bibliographystyle{siamplain}
\bibliography{Bibliography/Mine,Bibliography/NumericalMethodsHyperbolic,Bibliography/Books,Bibliography/ShallowWater,Bibliography/AnalysisHyperbolic}

\cleardoublepage

\begin{center}
\textbf{\uppercase{\large Supplemental Information}}
\end{center}
\vspace{5mm}

\section{Introduction}

This supplemental information is largely dedicated to the extension of the results presented in the main text to the case of general TVD slope limiters. The general results for the depth reconstruction substantially reduce the number of additional results required for any particular slope limiter to produce a self-monotone reconstruction. They are presented in the order they were deduced, outlining the steps we went through to prove our results, which we intend as a guide for anyone designing a well balanced scheme using a different reconstruction (\ie not minmod) as a base. In particular, the approach of convex combination and suppression is applicable to high order schemes, which introduces significant additional complexity to the derivation of a suitable $\gamma$, but otherwise the derivation will be similar.

The supplemental material is structured as follows. Fist we discuss the general properties of TVD slope limiters in \cref{sec:suppTVD}. We then mirror \cref{sec:reconH,sec:recon_q,sec:reconWH,sec:recon_phi} with \cref{sec:supp_reconH,sec:supprecon_q,sec:suppreconWH,sec:supprecon_phi}, presenting generalised results and further discussion. Finally in \cref{sec:suppother} we present some miscellaneous results.
\section{TVD schemes for scalar problems}\label{sec:suppTVD}

Our well-balanced, self-monotone reconstruction makes use of the $\minmod$ slope limiter \cref{eqn:minmod_recon}, which is a TVD reconstruction. We generalise our results to TVD slope limiters, and here we overview some classical results for TVD schemes, along with presenting the consequences of our monotonicity conditions. As in \cref{sec:monotone} we consider a scalar function $v$ satisfying a scalar conservation law \cref{eqn:FV_conservation} with flux $F = f(v)$ and source $\Psi = 0$.

The gradient of the piecewise linear reconstruction in each cell is given by
\begin{subequations}
\begin{equation}
	[v_x]_j \eqdef \sigma(v_{j-1},v_{j},v_{j+1};\{x\}).
\end{equation}
The slope limiter $\sigma$ is assumed to produce a reconstruction that is symmetric under reflection, thus
\begin{equation*}
	\sigma(v_{j-1},v_{j},v_{j+1}; \ldots x_{j-1/2},x_{j+1/2} \ldots) = -\sigma(v_{j+1},v_{j},v_{j-1}; \ldots -x_{j+1/2},-x_{j-1/2} \ldots).
\end{equation*}
\end{subequations}
We now relate the TVD property to the reconstruction.
\begin{definition}
	A reconstruction is called \emph{TVD} if, for any set of grid cell values, generalised MUSCL schemes \cite{ar_Osher_1985} with an E-flux \cite{ar_Osher_1984} cause the total variation to decrease or remain constant as time advances.
\end{definition}
We will not use this requirement directly, rather we will make use of the following result from \cite{ar_Tadmor_1988}.
\begin{lemma} \label{thm:TVD_Tadmor}
	If
	\begin{equation}	\label{eqn:TVD_Tadmor}
	\begin{aligned}
		&&
		\min(v_{j-1},v_{j}) &\leq v_{j-1/2}^+ \leq \max(v_{j-1},v_{j}),
		\\ \text{and}&&
		\min(v_{j},v_{j+1}) &\leq v_{j+1/2}^- \leq \max(v_{j},v_{j+1})
	\end{aligned}\end{equation}
	then the reconstruction is TVD.
\end{lemma}

From this we deduce the following. The arbitrary parameters $\alpha_{j+1/2}^{\pm}$ are included so that the bound can be tightened independently on each cell. With $\alpha_{j+1/2}^{\pm}=1$ this result is equivalent to the upper bound for flux limiters found in \cite{ar_Sweby_1984}.
\begin{lemma}	\label{thm:PWB_linearTVD}
	If
	\begin{equation} 	\label{eqn:PWB_linearTVD}
	\begin{aligned}
		&&
		\sigma &\geq \frac{2}{\Delta x_j} \max \pbrk*{\min\pbrk*{  \alpha_{j-1/2}^+\Delta v_{j-1/2} , 0 } , \min\pbrk*{  \alpha_{j+1/2}^-\Delta v_{j+1/2} , 0 }}
		\\	\text{and} &&
		\sigma &\leq \frac{2}{\Delta x_j} \min \pbrk*{\max\pbrk*{  \alpha_{j-1/2}^+\Delta v_{j-1/2} , 0 } , \max\pbrk*{  \alpha_{j+1/2}^-\Delta v_{j+1/2} , 0 }}
	\end{aligned}\end{equation}
	for some $\alpha_{j-1/2}^+,\alpha_{j+1/2}^-:\{x\} \rightarrow [0,1]$ then the slope limiter is TVD.
\end{lemma}
\paragraph{Remark} The reconstruction $\hat{v}$ satisfies
\begin{align}
	\abs{\hat{v}(x) - v_j} &\leq
	\min\ppar*{\alpha_{j-1/2}^+ \abs{\Delta v_{j-1/2}} , \alpha_{j+1/2}^- \abs{\Delta v_{j+1/2}} }
	&\text{for}&&
	x &\in C_j
\end{align}
That is, the constant $\alpha_{j+1/2}^-$ tells us how much the reconstruction is permitted to vary from the cell averaged value due to $\Delta v_{j+1/2}$, and similarly for $\alpha_{j-1/2}^+$ and $\Delta v_{j-1/2}$.
\paragraph{Remark} The $\minmod$ slope limiter \cref{eqn:minmod_recon} is most tightly bound when the $\alpha_{j+1/2}^\pm$ from \cref{thm:PWB_linearTVD} are equal to the $\alpha_{j+1/2}^\pm$ from \cref{eqn:minmod_recon}, hence this choice of notation.
\begin{proof}[Proof of \cref{thm:PWB_linearTVD}]
	We show the equivalence of \cref{eqn:PWB_linearTVD} with $\alpha_{j+1/2}^{\pm}=1$ to the requirements \cref{eqn:TVD_Tadmor} from \cref{thm:TVD_Tadmor} case by case. We first examine the various cases of the two requirements in \cref{eqn:TVD_Tadmor}.
	\begin{itemize}[left=0mm]
		\item The first requirement is $\min(v_{j-1},v_{j}) \leq v_{j-1/2}^+ \leq \max(v_{j-1},v_{j})$.
		\begin{subequations}\begin{flalign*}
			\textrm{When } v_{j-1} \leq v_{j}: &&
			v_{j-1} &\leq v_{j-1/2}^+  = v_j - \frac{\Delta x_j}{2} \sigma \leq v_j
			&\Leftrightarrow&&
			0 &\leq \sigma \leq \frac{2 \Delta v_{j-1/2}}{\Delta x_j}.
			\\
			\textrm{When } v_{j-1} \geq v_{j}: &&
			v_{j} &\leq v_{j-1/2}^+  = v_j - \frac{\Delta x_j}{2} \sigma \leq v_{j-1}
			&\Leftrightarrow&&
			%\frac{2 \Delta v_{j-1/2}}{\Delta x_j} &\leq \sigma \leq 0.&
			0 &\geq  \sigma \geq \frac{2 \Delta v_{j-1/2}}{\Delta x_j}.
		\end{flalign*}\end{subequations}
		\item The second requirement is $\min(v_{j},v_{j+1}) \leq v_{j+1/2}^- \leq \max(v_{j},v_{j+1})$.
		\begin{subequations}\begin{flalign*}
			\textrm{When } v_{j} \leq v_{j+1}: &&
			v_{j} &\leq v_{j+1/2}^-  = v_j + \frac{\Delta x_j}{2} \sigma \leq v_{j+1}
			&\Leftrightarrow&&
			0 &\leq \sigma \leq \frac{2 \Delta v_{j+1/2}}{\Delta x_j}.
			\\
			\textrm{When } v_{j} \geq v_{j+1}: &&
			v_{j+1} &\leq v_{j+1/2}^-  = v_j + \frac{\Delta x_j}{2} \sigma \leq v_j
			&\Leftrightarrow&&
			%\frac{2 \Delta v_{j+1/2}}{\Delta x_j} &\leq \sigma \leq 0.
			0 &\geq \sigma \geq \frac{2 \Delta v_{j+1/2}}{\Delta x_j}.
			\end{flalign*}\end{subequations}
	\end{itemize}
	We then combine these by considering all three of the cell values simultaneously
	\begin{itemize}[left=0mm]
		\item For $v_j \geq \max[v_{j-1},v_{j+1}]$ or $v_j \leq \min[v_{j-1},v_{j+1}]$ \cref{eqn:TVD_Tadmor,eqn:PWB_linearTVD} are both equivalent to $\sigma=0$.
		\item For $v_{j-1} \leq v_j \leq v_{j+1}$ \cref{eqn:TVD_Tadmor,eqn:PWB_linearTVD} are both equivalent to
		\begin{equation*}
			0 \leq \sigma \leq \min\pbrk*{ \frac{2 \Delta v_{j-1/2}}{\Delta x_j} , \frac{2 \Delta v_{j+1/2}}{\Delta x_j} }.
		\end{equation*}
		\item For $v_{j-1} \geq v_j \geq v_{j+1}$ \cref{eqn:TVD_Tadmor,eqn:PWB_linearTVD} are both equivalent to
		\begin{equation*}
			0 \geq \sigma \geq \max\pbrk*{ \frac{2 \Delta v_{j-1/2}}{\Delta x_j} , \frac{2 \Delta v_{j+1/2}}{\Delta x_j} }.
		\end{equation*}
	\end{itemize}
\end{proof}

Finally we present the consequences of imposing the monotonicity requirements (\cref{def:monotone}) on a slope limiter.

\begin{lemma}	\label{thm:PWB_selfmon}
	A slope limiter is self-monotone if and only if
	\begin{align}\label{eqn:PWB_sigma_selfmon}
		-\frac{2 \beta_j^-}{\Delta x_j} \leq \pdv{\sigma}{v_j} &\leq \frac{2 \beta_j^{+}}{\Delta x_j}
		&\textrm{for some}&&
		%\beta_j^{\pm}:(x_{j-3/2},x_{j-1/2},x_{j+1/2},x_{j+3/2}) \mapsto [0,1],
		\beta_j^{\pm}: \{x\} \mapsto [0,1],
	\end{align}
	similarly it is neighbour-monotone if and only if 
	\begin{align}\label{eqn:PWB_sigma_neighmon}
		\pdv{\sigma}{v_{j-1}} \leq 0&
		&\textrm{and}&&
		\pdv{\sigma}{v_{j+1}} \geq 0.
	\end{align}
\end{lemma}
\begin{proof}
	Self-monotonicity is defined as
	\begin{align}
		\pdv{v_{j-1/2}^{+}}{v_{j}} = 1 - \frac{\Delta x}{2} \pdv{\sigma}{v_j} &\geq 0,
		&\textrm{and}&&
		\pdv{v_{j+1/2}^{-}}{v_{j}} = 1 + \frac{\Delta x}{2} \pdv{\sigma}{v_j} &\geq 0
	\end{align}
	which are equivalent to \cref{eqn:PWB_sigma_selfmon} with $\beta_j^{\pm} = 1$.
	Neighbour-monotonicity is defined as
	\begin{align}
		\pdv{v_{j-1/2}^{+}}{v_{j-1}} = -\frac{\Delta x}{2} \pdv{\sigma}{v_{j-1}} &\geq 0,
		&\textrm{and} &&
		\pdv{v_{j+1/2}^{-}}{v_{j+1}} = \frac{\Delta x}{2} \pdv{\sigma}{v_{j+1}} &\geq 0
	\end{align}
	which are equivalent to \cref{eqn:PWB_sigma_neighmon}.
\end{proof}
\section{Reconstruction of depth}\label{sec:supp_reconH}

In this section we generalise the results of \cref{thm:results_h,sec:reconH}, guiding the reader through the reasoning used to deduce the structure of $\gamma_j$. We begin in \cref{sec:suppreconH_TVD} with a discussion of a general TVD slope limiter, and provide results for this generalisation. However, it has been found that to obtain a self-monotone reconstruction quite detailed information is required. For this reason we restrict ourselves to the specific case of the $\minmod$ slope limiter in \cref{sec:suppreconH_minmod}, where we derive the form of $\gamma_j$ included in \cref{thm:results_h}.
\subsection{Results for a general TVD slope limiter} \label{sec:suppreconH_TVD}

For the purpose of this subsection alone we utilise the generalised expressions
\begin{subequations}
\begin{align}\label{eqn:suppreconH_gradient}
	[h_x]_j^h 		&\eqdef \sigma_j^h
	&\textrm{where}&&
	\sigma_j^h &\eqdef \sigma^h(h_{j-1},h_{j},h_{j+1}; \{x\}),
	\\
	[h_x]_j^\eta 	&\eqdef \sigma_j^\eta - [b_x]_j
	&\textrm{where}&&
	\sigma_j^\eta &\eqdef \sigma^\eta(\eta_{j-1},\eta_{j},\eta_{j+1}; \{x\}),
\end{align}\end{subequations}
and $\sigma^h$ and $\sigma^\eta$ are symmetric, TVD, and CDP slope limiters, with coefficients for \cref{thm:PWB_linearTVD} taking least values $\alpha_{j+1/2}^{h\pm}$ and $\alpha_{j+1/2}^{\eta\pm}$ respectively, and coefficients for \cref{thm:PWB_selfmon} taking least values $\beta_{j}^{h\pm}$ and $\beta_{j}^{\eta\pm}$ respectively. The values of $h_j^\downarrow$ and $h_j^\uparrow$ are calculated using the parameters for the $\eta$ reconstruction, that is
\begin{equation}\begin{aligned}
	h_j^\downarrow &\eqdef \min \pbrk*{ h_j - \alpha_{j-1/2}^{\eta +}\Delta h_{j-1/2},h_j,h_j + \alpha_{j+1/2}^{\eta -}\Delta h_{j+1/2} }, 	\label{eqn:suppreconH_hmin}
	\\
	h_j^\uparrow &\eqdef \max \pbrk*{ h_j - \alpha_{j-1/2}^{\eta +}\Delta h_{j-1/2},h_j,h_j + \alpha_{j+1/2}^{\eta -}\Delta h_{j+1/2} }.
\end{aligned}\end{equation}

The first property we wish to establish is positivity. The reconstruction based on $h$ is always positive, whilst the one based on $\eta$ is able to become negative. Considering the depth $h_{j+1/2}^{\eta-}$, a lower bound can be constructed using \cref{thm:PWB_linearTVD},
\begin{equation}
	h_{j+1/2}^{\eta-} \geq h_j - \Delta b_j/2 + \alpha_{j+1/2}^{\eta-} \min(\Delta \eta_{j+1/2},0),
\end{equation}
Using that
\begin{subequations}\begin{align}
	h_j - \Delta b_j/2 + \alpha_{j+1/2}^{\eta-}\Delta\eta_{j+1/2} &\geq 0
	& \Leftrightarrow &&
	\frac{\Delta b_j/2 - \alpha_{j+1/2}^{\eta-}\Delta b_{j+1/2}}{h_j + \alpha_{j+1/2}^{\eta-}\Delta h_{j+1/2}} \leq 1,
	\\
	h_j - \Delta b_j/2 &\geq 0
	& \Leftrightarrow &&
	\frac{\Delta b_j/2}{h_j} \leq 1.
\end{align}\end{subequations}
and the symmetric property of the slope limiter we obtain
\begin{lemma}	\label{thm:suppreconH_bound_eta}	\begin{subequations}
	Let
	\begin{align*}
		\hat{\xi}_j &= \max \pbrk*{ \frac{-\Delta b_j/2 + \alpha_{j-1/2}^{\eta+}\Delta b_{j-1/2}}{h_j - \alpha_{j-1/2}^{\eta+}\Delta h_{j-1/2}} , \frac{\abs*{\Delta b_j}/2}{h_j} , \frac{\Delta b_j/2 - \alpha_{j+1/2}^{\eta-}\Delta b_{j+1/2}}{h_j + \alpha_{j+1/2}^{\eta-}\Delta h_{j+1/2}} },
		\\
		H_{j-1/2}^{+} &= h_j + \Delta b_j/2 - \alpha_{j-1/2}^{\eta+} \max[\Delta\eta_{j-1/2},0],
		\\
		H_{j+1/2}^{-} &= h_j - \Delta b_j/2 + \alpha_{j+1/2}^{\eta-} \min[\Delta\eta_{j+1/2},0],
	\end{align*}
	then [$H_{j-1/2}^{+} \geq 0$ and $H_{j+1/2}^{-} \geq 0$] if and only if $\hat{\xi}_j \leq 1$. Here $h_{j-1/2}^{\eta+} \geq H_{j-1/2}^{+}$ and $h_{j+1/2}^{\eta-} \geq H_{j+1/2}^{-}$ are bounds on the reconstruction.
\end{subequations}\end{lemma}
To aid analysis we define the local measure of depth relative to the bed variation
\begin{subequations}\label{eqn:suppreconH_xi}\begin{align}
	\xi_j &\eqdef \frac{h_j^\downarrow}{\Delta b_j^\uparrow}
\intertext{which takes the value $+\infty$ when $\Delta b_j^\uparrow=0$, where}
	%h_j^\downarrow &= \min\left[h_j - \alpha_{j-1/2}^{\eta+}\Delta h_{j-1/2},h_j,h_j + \alpha_{j+1/2}^{\eta-}\Delta h_{j+1/2}\right],	\label{eqn:suppreconH_hmin}
	%\\
	\Delta b_j^\uparrow &\geq \max\pbrk*{ \abs*{\Delta b_j/2-\!\alpha_{j-1/2}^{\eta+}\Delta b_{j-1/2}},\abs*{\Delta b_j/2},\abs*{\Delta b_j/2-\!\alpha_{j+1/2}^{\eta-}\Delta b_{j+1/2}} }. \label{eqn:suppreconH_dbmax}
\end{align}\end{subequations}
is (at least) the largest change in bed elevation across half the cell. We relate $\xi_j$ to $\hat{\xi}_j$ from \cref{thm:suppreconH_bound_eta} by
\begin{equation*}\begin{split}
	\hat{\xi}_j 
	&\leq \max \pbrk*{ \frac{\abs*{-\Delta b_j/2 + \alpha_{j-1/2}^{\eta+}\Delta b_{j-1/2}}}{h_j - \alpha_{j-1/2}^{\eta+}\Delta h_{j-1/2}} , \frac{\abs*{\Delta b_j}/2}{h_j} , \frac{\abs*{\Delta b_j/2 - \alpha_{j+1/2}^{\eta-}\Delta b_{j+1/2}}}{h_j + \alpha_{j+1/2}^{\eta-}\Delta h_{j+1/2}} }
	\\
	&\leq \frac{\max \pbrk*{ \abs*{ \Delta b_j/2-\alpha_{j-1/2}^{\eta+}\Delta b_{j-1/2}} , \abs{ \Delta b_j/2 } , \abs*{ \Delta b_j/2 - \alpha_{j+1/2}^{\eta-}\Delta b_{j+1/2} } } }
	{\min \pbrk*{ h_j - \alpha_{j-1/2}^{\eta+}\Delta h_{j-1/2},h_j,h_j + \alpha_{j+1/2}^{\eta-}\Delta h_{j+1/2} } }
	\leq \xi_j^{-1}.
\end{split}\end{equation*}

\begin{lemma}	\label{thm:suppreconH_def_xi}
	If $\xi_j \geq 1$ then $h_{j-1/2}^{\eta+} \geq 0$ and $h_{j+1/2}^{\eta-} \geq 0$ .
\end{lemma}

By \cref{thm:suppreconH_def_xi} the construction $\gamma_j = \gamma(\xi_j;\{x\})$ is reasonable, and we impose that $\Delta b_j^\uparrow$ is independent of $\{h\}$. In addition, we assume that $\gamma_j$ is invariant under reflection in $x$. Thus by \cref{thm:suppreconH_def_xi}, $\gamma_j$ is unrestricted by positivity constraints for $\xi_j \geq 1$. The expression for $\gamma_j$ we construct will be continuous so that $\gamma_j=0$ for $\xi_j \leq 1$, $\gamma_j=1$ for $\xi_j \geq \xi_j^{C}$, and non-decreasing within $1 \leq \xi_j \leq \xi_j^{C}$, where $\xi_j^{C}$ is a function from $\{x\}$ to $(1,\infty)$ and represents the critical value of $\xi_j$ above which the reconstruction is in $\eta$ only.

We turn our attention to the monotonicity properties, and by reflective symmetry we discuss only $h_{j+1/2}^-$. The conditions for self and neighbour-monotonicity are $\pdv*{h_{j+1/2}^{-}}{h_{j}} \geq 0$ and $\pdv*{h_{j+1/2}^{-}}{h_{j+1}} \geq 0$ respectively, which are equivalent to 
\begin{subequations}\begin{align}
	\pdv{\gamma_j}{\xi_j}  R_j &\leq S_j
	&\textrm{and}&&
	\pdv{\gamma_j}{\xi_j}  R_j &\leq N_j
	\label{eqn:suppreconH_selfneigh_inequal}
\end{align}
respectively. Here
\begin{align}
	R_j 
	&\eqdef \frac{h_{j+1/2}^{h-} - h_{j+1/2}^{\eta-}}{\Delta b_j^\uparrow} 
	= \frac{\Delta b_j/2 + (\Delta x_j/2)(\sigma_j^{h}-\sigma_j^{\eta})}{\Delta b_j^\uparrow}, 
	\label{eqn:suppreconH_defn_R}
	\\
	\begin{split}
		S_j
		&\eqdef \ppar*{\pdv{h_j^\downarrow}{h_j}}^{-1} \cdot \ppar*{ (1-\gamma_j) \pdv{h_{j+1/2}^{h-}}{h_j} + \gamma_j \pdv{h_{j+1/2}^{\eta-}}{h_j} }	\\
		&= \ppar*{\pdv{h_j^\downarrow}{h_j}}^{-1} \cdot \ppar*{ 1 + \frac{\Delta x_j}{2} \ppar*{  (1-\gamma_j) \pdv{\sigma_j^{h}}{h_j} + \gamma_j \pdv{\sigma_j^{\eta}}{\eta_j} } },
	\end{split}
	\label{eqn:suppreconH_defn_GammaS}
	\\
	\begin{split}
		N_j
		&\eqdef \ppar*{\pdv{h_j^\downarrow}{h_{j+1}} }^{-1} \cdot  \ppar*{ (1-\gamma_j) \pdv{h_{j+1/2}^{h-}}{h_{j+1}} + \gamma_j \pdv{h_{j+1/2}^{\eta-}}{h_{j+1}} }	\\
		&= \ppar*{\pdv{h_j^\downarrow}{h_{j+1}} }^{-1} \cdot \frac{\Delta x_j}{2} \ppar*{ (1-\gamma_j) \pdv{\sigma_j^{h}}{h_{j+1}} + \gamma_j \pdv{\sigma_j^{\eta}}{\eta_{j+1}} },
	\end{split}
	\label{eqn:suppreconH_defn_GammaN}
\end{align}\end{subequations}
and if $\pdv*{h_j^\downarrow}{h_j}=0$ then $S_j=+\infty$, if $\pdv*{h_j^\downarrow}{h_{j+1}}=0$ then $N_j=+\infty$.  We can simplify the inequalities \cref{eqn:suppreconH_selfneigh_inequal} by bounding $R_j$, $S_j$ and $N_j$ and using that $\gamma_j$ is non-decreasing in $\xi_j$.

\begin{lemma}	\label{thm:suppreconH_selfmon_general} \begin{subequations}
	Suppose $R_j \leq \hat{R}_j$ is the least upper bound and $S_j\geq \hat{S}_j$, $N_j\geq \hat{N}_j$ are the greatest lower bounds, where each bound is a function of $(\gamma_j;\{x\})$. If 
	\begin{equation}	\label{thm:suppreconH_selfmon_general_S}
		\pdv{\gamma_j}{\xi_j} \leq \frac{\hat{S}_j}{\hat{R}_j},
	\end{equation}
	then the reconstruction is self-monotone, and if for some discretized values $R_j = \hat{R}_j$ and $S_j= \hat{S}_j$ simultaneously then \cref{thm:suppreconH_selfmon_general_S} is necessary for self-monotonicity. If 
	\begin{equation}	\label{thm:suppreconH_selfmon_general_N}
		\pdv{\gamma_j}{\xi_j} \leq \frac{\hat{N}_j}{\hat{R}_j},
	\end{equation}
	then the reconstruction is neighbour-monotone, and if for some discretized values $R_j = \hat{R}_j$ and $N_j= \hat{N}_j$ simultaneously then \cref{thm:suppreconH_selfmon_general_N} is necessary for neighbour-monotonicity.
\end{subequations}\end{lemma}
\Cref{thm:suppreconH_selfmon_general} is implied by the following result.
\begin{lemma}\label{thm:suppreconH_inequalities}
	Suppose that for some $x \geq 0$, $y\in\mathscr{Y}$ we have functions $\Gamma: \mathscr{Y} \rightarrow \mathbb{R}_0^+$, $R: \mathscr{Y} \rightarrow \mathbb{R}$. The statement
	\begin{align*}
		x R(y) &\leq \Gamma(y) \quad \forall \quad y \in \mathscr{Y}
		&&\textrm{is equivalent to}&
		x &\leq \inf_{y:R(y)>0} \frac**{\Gamma(y)}{R(y)},
	\end{align*}
	moreover, if the minimal value of $\Gamma$ and maximal value of $R$ is attained at the same $y$, then both are equivalent to 
	\begin{equation*}
		x \leq \frac**{\inf_{y:R(y)>0}(\Gamma(y))}{\sup_{y:R(y)>0}(R(y))}.
	\end{equation*}
\end{lemma}
To obtain an explicit expression for $\gamma_j$ we construct bounds for $R_j$, $S_j$, and $N_j$.

\begin{lemma}	\label{thm:suppreconH_bound_R_general}
	We neglect the cases [$h_{j-1} < h_{j} < h_{j+1}$ and $\eta_{j-1} < \eta_{j} < \eta_{j+1}$] and [$h_{j-1} > h_{j} > h_{j+1}$ and $\eta_{j-1} > \eta_{j} > \eta_{j+1}$] and assume $\alpha_{j-1/2}^{\eta+}=\alpha_{j-1/2}^{h+}$ and $\alpha_{j+1/2}^{\eta-}=\alpha_{j+1/2}^{h-}$. $R_j\leq 1$ is an upper bound for $R_j$, and provided $\Delta b_j^\uparrow=\Delta b_j/2$ is possible for any $(\gamma_j;\{x\})$ then this is the least upper bound.
\end{lemma}
\begin{proof}\begin{subequations}
	We start by showing that $R_j$ can equal the bound. When both $h$ and $\eta$ are at an extrema, that is $\Delta h_{j-1/2}\Delta h_{j+1/2} \leq 0$ and $\Delta\eta_{j-1/2} \Delta\eta_{j+1/2} \leq 0$, we have $R_j = \Delta b_j / 2 \Delta b_{j}^\uparrow$. Clearly $R_j \leq 1$ by the definition of $\Delta b_{j}^\uparrow$, and when $\Delta b_{j}^\uparrow = \Delta b_j / 2$ we have $R = 1$. We now show that this is an upper bound for the other cases included in the lemma.
	\begin{itemize}[left=0mm]
		\item If [$\Delta h_{j-1/2} \leq 0$ or $\Delta h_{j+1/2} \leq 0$] and [$\Delta \eta_{j-1/2} \geq 0$ or $\Delta \eta_{j+1/2} \geq 0$] then
		\begin{equation}
		\eqnwide
			R_j
			\leq \frac{\Delta b_j/2}{\Delta b_j^\uparrow}
			\leq 1.
		\end{equation}
		\item If $\Delta h_{j+1/2} \geq 0$ and $\Delta\eta_{j+1/2} \leq 0$ then
		\begin{equation}\label{eqn:thm_suppreconH_bound_R_boundL}
		\eqnwide
			R_j
			\leq \frac{\Delta b_j/2 + \alpha_{j+1/2}^{h-} (\Delta h_{j+1/2} - \Delta\eta_{j+1/2})}{\Delta b_j^\uparrow}
			= \frac{\Delta b_j/2 - \alpha_{j+1/2}^{h-} \Delta b_{j+1/2}}{\Delta b_j^\uparrow}
			\leq 1.
		\end{equation}
		\item If $\Delta h_{j-1/2} \geq 0$ and $\Delta\eta_{j-1/2} \leq 0$ then
		\begin{equation}\label{eqn:thm_suppreconH_bound_R_boundR}
		\eqnwide
			R_j
			\leq \frac{\Delta b_j/2 + \alpha_{j-1/2}^{h+} (\Delta h_{j-1/2} - \Delta\eta_{j-1/2})}{\Delta b_j^\uparrow}
			= \frac{\Delta b_j/2 - \alpha_{j-1/2}^{h+} \Delta b_{j-1/2}}{\Delta b_j^\uparrow}
			\leq 1.
		\end{equation}
	\end{itemize}
\end{subequations}\end{proof}

\begin{lemma}	\label{thm:suppreconH_bound_S}
	$S_j \geq \hat{S}_j$ is a lower bound for $S_j$, where
	\begin{equation*}
		\hat{S}_j(\gamma_j) =
		\min \pbrk*{ 1 - \gamma_j\beta_j^{\eta-} , \frac{1 - (1-\gamma_j) \beta_j^{h-} - \gamma_j \beta_j^{\eta-}}{1-\alpha_{j-1/2}^{\eta+}} , \frac{ 1 - (1-\gamma_j) \beta_j^{h-} - \gamma_j \beta_j^{\eta-}}{1-\alpha_{j+1/2}^{\eta-}} }.
	\end{equation*}
	Note that in this expression we treat $1/(1-\alpha) = +\infty$ for $\alpha = 1$.
	Moreover, if for any $(\gamma_j;\{x\})$ there are some conditions under which $h_j^\downarrow = h_j$ and $\sigma^{\eta}$ satisfies the lower bound in \cref{eqn:PWB_sigma_selfmon} as equality, and there are some other conditions under which $h_j^\downarrow$ takes the other expressions in the $\min$ \cref{eqn:suppreconH_hmin} and both $\sigma^{h}$ and $\sigma^{\eta}$ equal the lower bound, then this is the greatest lower bound.
\end{lemma}
\begin{proof}
	We prove by considering the different values possible for $h_j^\downarrow$
	\begin{itemize}[left=0mm]
		\item If $h_j^\downarrow = h_j$ then $h_j \leq \min[h_{j-1},h_{j+1}]$ and $\sigma_j^{h} = 0$, thus
		\begin{equation*}
			S_j \geq 1 - \gamma_j\beta_j^{\eta-}
		\end{equation*}
		and this is an equality when $\pdv*{\sigma_j^{\eta}}{\eta_j} = -2 \beta_j^{\eta-}/\Delta x_j$, thus is a greatest lower bound for this case.
		
		\item If $h_j^\downarrow = h_j - \alpha_{j-1/2}^{\eta+}\Delta h_{j-1/2}$ then, so long as $\alpha_{j-1/2}^{\eta+} \neq 1$,
		\begin{equation*}
			S_j \geq \frac{1 - (1-\gamma_j) \beta_j^{h-} - \gamma_j \beta_j^{\eta-}}{1-\alpha_{j-1/2}^{\eta+}}
		\end{equation*}
		and this is an equality when $\pdv*{\sigma_j^{h}}{h_j} = -2 \beta_j^{h-}/\Delta x_j$ and  $\pdv*{\sigma_j^{\eta}}{\eta_j}= -2 \beta_j^{\eta-}/\Delta x_j$, thus is a greatest lower bound for this case.
		
		\item The $h_j^\downarrow = h_j + \alpha_{j+1/2}^{\eta-}\Delta h_{j+1/2}$ case is equivalent to the above by symmetry, and yields the same bound with $\alpha_{j+1/2}^{\eta-}$ substituted for $\alpha_{j-1/2}^{\eta+}$.
	\end{itemize}
\end{proof}

\begin{lemma}	\label{thm:suppreconH_bound_N}
	If, for any $(\gamma_j;\{x\})$, $\alpha_{j+1/2}^{\eta-} \neq 0$ and there is some condition under which $h_j^\downarrow = h_j + \alpha_{j+1/2}^{\eta-}\Delta h_{j+1/2}$ and $\sigma_j^{h}$ and $\sigma_j^{\eta}$ simultaneously satisfy the lower bound in \cref{eqn:PWB_sigma_selfmon} as equalities, then $N_j \geq \hat{N}_j = 0$ is the greatest lower bound.
\end{lemma}
\begin{proof}
	The only condition under which $\pdv*{h_j^\downarrow}{h_{j+1}} \neq 0$ is when $\alpha_{j+1/2}^{\eta-} \neq 0$ and $h_j^\downarrow = h_j + \alpha_{j+1/2}^{\eta-}\Delta h_{j+1/2}$, in which case $N_j\geq 0$ and this an equality when $\pdv*{\sigma_j^{h}}{h_{j+1}} = 0$ and  $\pdv*{\sigma_j^{\eta}}{\eta_{j+1}}= 0$.
\end{proof}

\Cref{thm:suppreconH_bound_R_general} omits the case of $h$ and $\eta$ varying monotonically in the same direction, thus the bound is not general. It has not been found possible to extend the methods used in \cref{thm:suppreconH_bound_R_general} to this case for general TVD slope limiters. \Cref{thm:suppreconH_bound_N} indicates that, for a wide class of slope limiters, $\gamma_j$ cannot be both an increasing function of $\xi_j$ and neighbour-monotone. This prevents the transition from reconstruction in $h$ to reconstruction in $\eta$. We see that we will have to drop the neighbour-monotone condition to obtain a well-balanced reconstruction.
\subsection{Results for the minmod slope limiter} \label{sec:suppreconH_minmod}

We continue by employing a minmod slope limiter \cref{eqn:minmod_recon}. One may desire to use different parameters for the two reconstructions: $\alpha_j^h$ and $\alpha_{j+1/2}^{h\pm}$ for the reconstruction in $h$; and $\alpha_j^\eta$ and $\alpha_{j+1/2}^{\eta\pm}$ for the reconstruction in $\eta$, thus $\beta_j^{\omitdummy-} = \alpha_{j+1/2}^{\omitdummy-}$ and $\beta_j^{\omitdummy+} = \alpha_{j-1/2}^{\omitdummy+}$ where $\omitdummy$ stands for $h$ or $\eta$. We begin by showing why the parameters must be the same.

\begin{lemma}	\label{thm:suppreconH_nobound_R_minmod}
	Suppose $\alpha_{j-1/2}^{\omitdummy+} , \alpha_j^{\omitdummy} , \alpha_{j+1/2}^{\omitdummy-} > 0$. If $\alpha_{j-1/2}^{h+} \neq \alpha_{j-1/2}^{\eta+}$ or $\alpha_{j+1/2}^{h-} \neq \alpha_{j+1/2}^{\eta-}$ then $R_j$ has no upper bound of the form required for \cref{thm:suppreconH_selfmon_general}, neither does it if $(\alpha_{j-1/2}^{\omitdummy+})^{-1} + (\alpha_{j+1/2}^{\omitdummy-})^{-1} < (\alpha_{j}^{\omitdummy})^{-1}$ and  $\alpha_{j}^{h} \neq \alpha_{j}^{\eta}$.
\end{lemma}
\begin{proof}
	We first consider when both $h$ and $\eta$ are reconstructed using left differences
	\begin{gather*}
		\begin{aligned}
			\sigma_j^{h} &= \frac{2 \alpha_{j-1/2}^{h+} \Delta h_{j-1/2}}{\Delta x_{j}}
			& \textrm{and} &&
			\sigma_j^{\eta} &= \frac{2 \alpha_{j-1/2}^{\eta+} \Delta \eta_{j-1/2}}{\Delta x_{j}},
		\end{aligned}
		\\
		\begin{aligned}
			&\text{thus}& 
			R_j &= \frac{\Delta b_j/2 - \alpha_{j-1/2}^{\eta+} \Delta b_{j-1/2} + (\alpha_{j-1/2}^{h+}-\alpha_{j-1/2}^{\eta+})\Delta h_{j-1/2}}{\Delta b_j^\uparrow}.
		\end{aligned}
	\end{gather*}
	The value of $\Delta h_{j-1/2}$ cannot be upper or lower bounded by a function of $\gamma_j$ and $\{x\}$. Therefore, if $\alpha_{j-1/2}^{h+} \neq \alpha_{j-1/2}^{\eta+}$ then $R_j$ has no upper bound.
	
	Performing equivalent analysis for when both $h$ and $\eta$ are reconstructed using right differences obtains the result that, if $\alpha_{j+1/2}^{h-} \neq \alpha_{j+1/2}^{\eta-}$ then $R_j$ has no upper bound. For the case of centred differences the condition $(\alpha_{j-1/2}^{\omitdummy+})^{-1} + (\alpha_{j+1/2}^{\omitdummy-})^{-1} < (\alpha_{j}^{\omitdummy})^{-1}$ ensures that they will be used (\cref{thm:minmod_centdiff_condn}). The rest of the proof is equivalent to that of left differences.
\end{proof}

Because of \cref{thm:suppreconH_nobound_R_minmod} we take $\alpha_j^\eta = \alpha_j^h$ and $\alpha_{j+1/2}^{\eta\pm} = \alpha_{j+1/2}^{h\pm}$ for the remainder of this subsection. This makes it possible to construct general bounds. In addition, we will assume that 
\begin{equation}\label{eqn:suppreconH_dbmax_strong}
\begin{split}
	\Delta b_j^\uparrow &\geq \max\left[|\Delta b_j/2-\alpha_{j-1/2}^{h+}\Delta b_{j-1/2}|,|\Delta b_j/2|,
	\right. \\ & \qquad \left.
	|\Delta b_j/2 - \alpha_j^h (b_{j+1}-b_{j-1})|,|\Delta b_j/2 - \alpha_{j+1/2}^{h-}\Delta b_{j+1/2}|\right]. 
\end{split}\end{equation}
(the third expression in the $\max$ was not included in \cref{eqn:suppreconH_dbmax}) which enables the following result.

\begin{lemma}	\label{thm:suppreconH_bound_R_minmod}
	The variable $R_j$ has bounds $-1 \leq R_j \leq 1$. Provided that, for any $(\gamma_j;\{x\})$, $\Delta b_j^\uparrow=\Delta b_j/2$ is possible then this is the leat upper bound.
\end{lemma}
\begin{proof}
	The fact that $R=1$ is possible comes from \cref{thm:suppreconH_bound_R_general}, along with the bound for all but a few cases which we consider here. We consider first when $h_{j-1} < h_{j} < h_{j+1}$ and $\eta_{j-1} < \eta_{j} < \eta_{j+1}$.
	\begin{itemize}[left=0mm]
		\item If $\sigma_j^{\eta}\Delta x_j/2 = \alpha_{j-1/2}^{h+} \Delta \eta_{j-1/2}$ then by $\sigma_j^{h}\Delta x_j/2 \leq \alpha_{j-1/2}^{h+} \Delta h_{j-1/2}$ we obtain \cref{eqn:thm_suppreconH_bound_R_boundL}.
		
		\item If $\sigma_j^{\eta}\Delta x_j/2 = \alpha_{j}^{h} (\eta_{j+1}-\eta_{j-1})$ then by $\sigma_j^{h}\Delta x_j/2 \leq \alpha_{j}^{h} (h_{j+1}-h_{j-1})$
		\begin{equation} \label{eqn:thm_suppreconH_bound_R_boundC}
		\eqnwide
			R_j
			\leq \frac{\Delta b_j/2 + \alpha_{j}^{h} \ppar*{h_{j+1}-h_{j-1}-\eta_{j+1}+\eta_{j-1}}}{\Delta b_j^\uparrow}
			= \frac{\Delta b_j/2 - \alpha_{j}^{h} \ppar*{b_{j+1} - b_{j-1}}}{\Delta b_j^\uparrow}
			\leq 1.
		\end{equation}
		
		\item If $\sigma_j^{\eta}\Delta x_j/2 = \alpha_{j+1/2}^{h-} \Delta \eta_{j+1/2}$ then by $\sigma_j^{h}\Delta x_j/2 \leq \alpha_{j+1/2}^{h-} \Delta h_{j+1/2}$ we obtain \cref{eqn:thm_suppreconH_bound_R_boundR}
	\end{itemize}
	Consider now $h_{j-1} > h_{j} > h_{j+1}$ and $\eta_{j-1} > \eta_{j} > \eta_{j+1}$.
	\begin{itemize}[left=0mm]
		\item If $\sigma_j^{h}\Delta x_j/2 = \alpha_{j-1/2}^{h+} \Delta h_{j-1/2}$ then by $\sigma_j^{\eta}\Delta x_j/2 \geq \alpha_{j-1/2}^{h+} \Delta \eta_{j-1/2}$ we obtain \cref{eqn:thm_suppreconH_bound_R_boundL}.
		
		\item If $\sigma_j^{h}\Delta x_j/2 =  \alpha_{j}^{h} (h_{j+1}-h_{j-1})$ then by $\sigma_j^{\eta}\Delta x_j/2 \geq \alpha_{j}^{h} (\eta_{j+1}-\eta_{j-1})$ we obtain \cref{eqn:thm_suppreconH_bound_R_boundC}.
		
		\item If $\sigma_j^{h}\Delta x_j/2 = \alpha_{j+1/2}^{h-} \Delta h_{j+1/2}$ then by $\sigma_j^{\eta}\Delta x_j/2 \geq \alpha_{j+1/2}^{h-} \Delta \eta_{j+1/2}$ we obtain \cref{eqn:thm_suppreconH_bound_R_boundR}.
	\end{itemize}
	The lower bound comes from symmetry. Under reflection $\sigma_j^{\omitdummy} \mapsto -\sigma_j^{\omitdummy}$, $\Delta b_j \mapsto -\Delta b_j$, and $\Delta b_j^\uparrow \mapsto \Delta b_j^\uparrow$, thus $R_j \mapsto -R_j$. Thus $-R_j \leq 1$ is a (least) upper bound, from which we deduce $R_j \geq -1$ is a (greatest) lower bound.
\end{proof}

\begin{lemma}		\label{thm:suppreconH_gamma_inequality_optimal} \begin{subequations}
	If
	\begin{equation}\label{eqn:suppreconH_gamma_inequality_optimal_self}
		\pdv{\gamma_j}{\xi_j} \leq \min\left[ 1 - \gamma_j \alpha_{j+1/2}^{h-} , \frac{1 - \alpha_{j+1/2}^{h-}}{1-\alpha_{j-1/2}^{h+}} \right].
	\end{equation}
	then the reconstruction is self-monotone, and if 
	\begin{equation}\label{eqn:suppreconH_gamma_inequality_optimal_neigh}
		\pdv{\gamma_j}{\xi_j} = 0.
	\end{equation}
	the reconstruction is neighbour-monotone. If discretized values of the fields other than $h$ and $b$ exist such that \cref{eqn:suppreconH_dbmax_strong} is an equality for all $\{h\}$ then \cref{eqn:suppreconH_gamma_inequality_optimal_self,eqn:suppreconH_gamma_inequality_optimal_neigh} are necessary.
\end{subequations}\end{lemma}
\begin{proof}
	Sufficiency is the result of \cref{thm:suppreconH_selfmon_general} using \cref{thm:suppreconH_bound_R_general,thm:suppreconH_bound_S,thm:suppreconH_bound_N,thm:suppreconH_bound_R_minmod}. Necessity requires that equality with the bounds for $S_j$ and $R_j$ is possible simultaneously, as well for $N_j$ and $R_j$. To show for $S_j$ and $R_j$ we consider the case $h_j \leq \min[h_{j-1},h_{j+1}]$, $\sigma_j^{\eta} = 2\alpha_{j+1/2}^{h-}\Delta \eta_{j+1/2}/\Delta x_j$, $\Delta b_j^\uparrow = \Delta b_j/2 - \alpha_{j+1/2}^{h-}\Delta b_{j+1/2}$, for which
	\begin{align}
		R_j &= 1
		&\textrm{and}&&
		S_j &= 1 - \gamma_j \alpha_{j+1/2}^{h-}.
	\end{align}
	We also consider the case $\sigma_j^{h} = 2\alpha_{j+1/2}^{h-}\Delta h_{j+1/2}/\Delta x_j$, $\sigma_j^{\eta} = 2\alpha_{j+1/2}^{h-}\Delta \eta_{j+1/2}/\Delta x_j$, $h_j^\downarrow = h_j - \alpha_{j-1/2}^{h+}\Delta h_{j-1/2}$ and $\Delta b_j^\uparrow = \Delta b_j/2 - \alpha_{j+1/2}^{h-}\Delta b_{j+1/2}$, for which
	\begin{align}
		R_j &= 1
		&\textrm{and}&&
		S_j &= \frac{1 - \alpha_{j+1/2}^{h-}}{1-\alpha_{j-1/2}^{h+}}.
	\end{align}
	To show for $N_j$ and $R_j$ we consider the case $\sigma_j^{h} = 2\alpha_{j-1/2}^{h+}\Delta h_{j-1/2}/\Delta x_j$, $\sigma_j^{\eta} = 2\alpha_{j-1/2}^{h+}\Delta \eta_{j-1/2}/\Delta x_j$, $h_j^\downarrow = h_j + \alpha_{j+1/2}^{h-}\Delta h_{j+1/2}$ and $\Delta b_j^\uparrow = \Delta b_j/2 - \alpha_{j-1/2}^{h+}\Delta b_{j-1/2}$, for which
	\begin{align}
		R_j &= 1
		&\textrm{and}&&
		N_j &= 0.
	\end{align}
\end{proof}

We observe that we cannot transition between the two reconstructions and be neighbour-monotone (at least for this parametrisation of $\gamma_j$), which is why we focus on self-monotonicity. By taking \cref{eqn:suppreconH_gamma_inequality_optimal_self} as equality we can produce an expression for $\gamma_j$, and thereby obtain a self-monotone reconstruction. However, our final goal here is to bound $u$ by producing a lower bound on the reconstructed $h$. To make this process simpler we assume that $0 < \alpha_{j+1/2}^{h\pm}<1$ and set
\begin{align} \label{eqn:PWB_gamma_inequality_simplified}
	\pdv{\gamma_j}{\xi_j} &= G_j \leq 1-\alpha_j^{h\uparrow}
	&\textrm{where}&&
	\alpha_j^{h\uparrow} &= \max[ \alpha_{j-1/2}^{h+} , \alpha_{j+1/2}^{h-} ]
\end{align}
over the region $1 < \xi_j < \xi_j^{C}$ where $\gamma_j$ is not constant. Note that $G_j$ is a function of $\{x\}$, and the inclusion of $\alpha_{j-1/2}^{h+}$ in $\alpha_j^{h\uparrow}$ is to provide the symmetry under reflection that is required. This results in the expression
\begin{align} \label{eqn:PWB_gamma_solution}
	\gamma_j(\xi_j) &=
	\begin{cases}
		0 & \textrm{if } \xi_j \leq 1	\\
		G_j(\xi_j - 1)	& \textrm{if } 1 \leq \xi_j \leq \xi_j^{C}	\\
		1 & \textrm{if } \xi_j^{C} \leq \xi_j
	\end{cases}
	&\textrm{where}&&
	\xi_j^{C} &= 1 + \frac{1}{G_j}
\end{align}
used in \cref{thm:results_h}.
\section{Reconstruction of flux} \label{sec:supprecon_q}

In this section we generalise the results of \cref{thm:results_u,sec:recon_q}. The flux $q$ is reconstructed as 
\begin{align*}
	[q_x]_j &= \kappa_j \sigma^{q}(q_{j-1},q_{j},q_{j+1};\{x\}) 
	&&\textrm{where}&
	\kappa_j &\eqdef \min \pbrk*{ 1 ,  \frac{K_{j-1/2}^+ h_j}{h_{j-1}} , \frac{K_{j+1/2}^- h_j}{h_{j+1}} },
\end{align*}
where $\sigma^{q}$ is a symmetric, TVD, and CDP reconstruction function bounded as in \cref{thm:PWB_linearTVD} with the coefficients having least values $\alpha_{j+1/2}^{q\pm}$. To produce a bound on $q_{j+1/2}^-$ ($q_{j-1/2}^+$ similar by symmetry) we first observe that
\begin{align*}\begin{split}
	\abs*{q_{j+1/2}^-} 
	&\leq \max \pbrk*{ \abs*{q_j} , \abs*{q_j + \kappa_j  \alpha_{j+1/2}^{q-} (q_{j+1}-q_j)}  }
	\\
	&\leq \max \pbrk*{ \abs*{u_j} h_j , (1 - \kappa_j \alpha_{j+1/2}^{q-}) \abs*{u_j} h_j + \kappa_j \alpha_{j+1/2}^{q-} \abs*{u_{j+1}} h_{j+1}}
\end{split}\end{align*}
Next we use that, no matter the discretized values of depth, $\kappa_j h_{j+1} \leq K_{j+1/2}^- h_j$, thus
\begin{align*}\begin{split}
	\abs*{q_{j+1/2}^-} 
	&\leq h_j \max \pbrk*{ \abs*{u_j} , (1 - \kappa_j \alpha_{j+1/2}^{q-}) \abs*{u_j} + K_{j+1/2}^- \alpha_{j+1/2}^{q-} \abs*{u_{j+1}}}
	\\
	&\leq h_j \ppar*{ \abs*{u_j} + K_{j+1/2}^- \alpha_{j+1/2}^{q-} \abs*{u_{j+1}} }.
\end{split}\end{align*}
To turn this into a bound on $\abs*{u_{j+1/2}^-}$ we make use of the bounds on the reconstructed depth from \cref{thm:reconH_h_bound}
\begin{align*}
	\frac{h_j}{h_{j+1/2}^-} \leq (G_j+1) \frac{h_j}{h_j^\downarrow} \leq \frac{G_j + 1}{1 - \alpha_j^{h \uparrow}}
\end{align*}
where the second inequality is from consideration of the case $h_{j-1} = h_{j+1} = 0$, thus
\begin{align}
	\abs*{u_{j+1/2}^-} \leq \frac{G_j + 1}{1 - \alpha_j^{h \uparrow}} \ppar*{ \abs*{u_j} + K_{j+1/2}^- \alpha_{j+1/2}^{q-} \abs*{u_{j+1}} }.
\end{align}
\section{Modified depth reconstruction for the inclusion of width} \label{sec:suppreconWH}

In this section we expand on the discussion presented in \cref{sec:reconWH}; specifically the construction of a positivity preserving reconstruction from a product of two piecewise-liner reconstructions.

\subsection{Quadratic and linear reconstructions} \label{sec:suppreconWH_simp}

We consider a simplified problem on $-1 \leq x \leq 1$ where we have two functions $f_1(x) \eqdef 1 + a_1 x$, $f_2(x) \eqdef 1 + a_2 x$, with $a_1,a_2 \in [-1,1]$. These two functions have a product $f_1(x) f_2(x) = 1 + (a_1+a_2)x + a_1 a_2 x^2$, the average value of which is not $1$, and so cannot be used a reconstruction. If a quadratic reconstruction is desired, then the average value of the product may be adjusted to produce the quadratic reconstruction
\begin{equation}\begin{split}
	f_q(x) &\eqdef f_1(x) f_2(x) - \int_{-1}^{1} f_1(x') f_2(x') \dd{x'} + 1 	\\
	&= 1 - \textfrac{1}{3}a_1 a_2 + (a_1 + a_2)x + a_1 a_2 x^2.
\end{split}\end{equation}
This function satisfies
\begin{align}
	f_q(-1) &= (1-a_1)(1-a_2) - \textfrac{1}{3}a_1 a_2,
	&
	f_q(1) &= (1+a_1)(1+a_2) - \textfrac{1}{3}a_1 a_2.
\end{align}
This if $\abs{a_1},\abs{a_2} \leq \alpha$ (\cf \cref{thm:PWB_linearTVD}) then, for positivity, we require
\begin{align}
	(1-\alpha)^2 - \textfrac{1}{3} \alpha^2 &\geq 0,
	&\text{thus}&&
	\alpha &\leq \frac{3 - \sqrt{3}}{2} = 0.634 \text{ (3 \sf)}.
\end{align}
This is quite a harsh restriction, seeing as $\alpha > 1/2$ is required for a second order convergence on a uniform grid and larger values will be required on a non-uniform grid, and indicates that the quadratic reconstruction is not the best approach.

\begin{figure}[t]
	\begin{center}
		\includegraphics{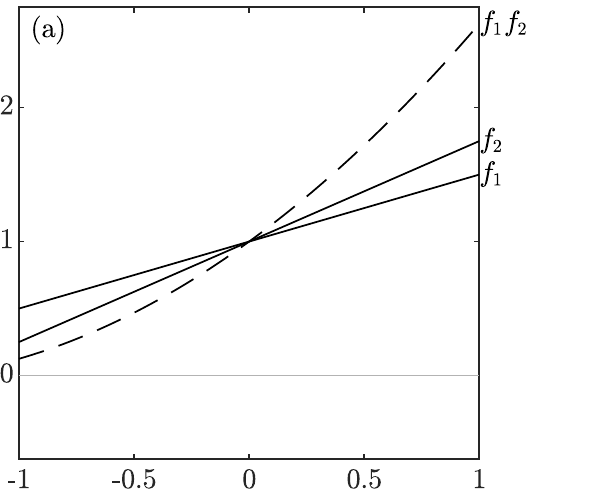}
		\includegraphics{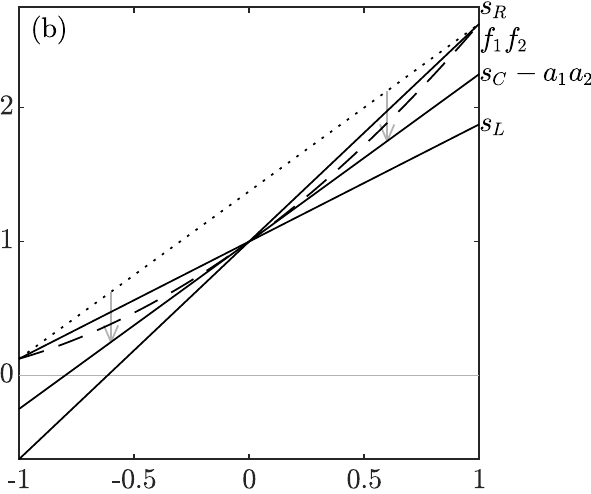}
	\end{center}
	\caption{Plots of the product reconstruction for our simplified example. In both figures we show the product of the functions $f_1(x)f_2(x)$ in dashed lines. In (a) we show the individual functions $f_1(x) = 1+x/2$ and $f_2(x)=1+3x/4$ as solid lines. In (b) we show the reconstructions $s_{L}(x)$, $s_{C}(x)-a_1 a_2$, and $s_{R}(x)$ as solid lines, along with the secant $s_{C}(x)$ as a dotted line.}
	\label{fig:suppreconWH_Secant}
\end{figure}

A linear reconstruction can be performed by computing some gradient from the product $f_1 f_2$, which may be any value in the range
\begin{equation}
	a_1 + a_2 - 2 \abs{a_1 a_2} \leq \dv{}{x} (f_1(x) f_2(x)) = a_1 + a_2 + 2 a_1 a_2 x \leq a_1 + a_2 + 2 \abs{a_1 a_2},
\end{equation}
thus the linear reconstruction is
\begin{equation}
	f_l (x) = 1 + \ppar*{ a_1 + a_2 + r a_1 a_2 } x
\end{equation}
for some $r \in [-2,2]$. The advantage of this approach is that it gives a tunable parameter $r$ which may be used to achieve a positive reconstruction. The choice $r=0$ yields an expression with the appearance of the product rule (\cf \cref{eqn:reconWH_grad_stucture}), indeed it is the gradient of $f_1 f_2$ at $x=0$, and thus may be an appealing option. The same gradient may also be achieved by the secant $s_C(x)$ though $f_1 f_2$ at $x=-1,1$, thus
\begin{align}
	\eval*{f_l(x)}_{r=0} &= s_C(x) - a_1 a_2
	&\text{where}&&
	s_C(x) &= 1 + a_1 a_2 + (a_1 + a_2) x
\end{align}
as depicted in \cref{fig:suppreconWH_Secant}. For $a_1 a_2 > 0$ the reconstruction is below the (positive) secant line, and thus may be negative, indeed for positivity we require $\abs{a_1 + a_2} \leq 1$. For the case $a_1 \geq 0$ and $a_2 \geq 0$ there is another secant that may be employed, the left scant through $f_1 f_2$ at $x=-1,0$,
\begin{align}
	s_L(x) &= 1 + (a_1 + a_2 - a_1 a_2) x.
\end{align}
This secant is the reconstruction with $r=-1$, and is positive for $1 + a_1 + b_1 - a_1 b_1$, thus for $a_1 \geq 0$ and $a_2 \geq 0$. A similar property holds for $a_1 \leq 0$ and $a_2 \leq 0$ and the right secant through $f_1 f_2$ at $x=0,1$, 
\begin{align}
	s_R(x) &= 1 + (a_1 + a_2 + a_1 a_2) x,
\end{align}
which is the reconstruction with $r=1$. In \cref{eqn:reconWH_grad} we have chosen $r = - \sign(a_1)$ so that the left secant is used when $f_1$ is increasing and the right when $f_1$ is decreasing. When $a_1=0$ 

\subsection{Bounds for the reconstruction} \label{sec:suppreconWH_bound}

Here we explicitly prove the bounds in \cref{eqn:reconWH_bounds}. First for \cref{eqn:reconWH_value} we use
\begin{subequations}\begin{gather}
	v_{j+1/2}^- = v_j + \frac{w_j^\downarrow}{w_{j+1/2}^-} \frac{\Delta x_j}{2} [v_x]_j
	\qquad
	\text{and}
	\qquad
	v_{j+1/2}^{v-} = v_j + \frac{\Delta x_j}{2} [v_x]_j,
\shortintertext{if $[v_x]_j \geq 0$ then}
	v_j 
	\leq v_j + \frac{w_j^\downarrow}{w_{j+1/2}^-} \frac{\Delta x_j}{2} [v_x]_j 
	\leq v_j + \frac{\Delta x_j}{2} [v_x]_j,
\shortintertext{and if $[v_x]_j \leq 0$ then}
	v_j 
	\geq v_j + \frac{w_j^\downarrow}{w_{j+1/2}^-} \frac{\Delta x_j}{2} [v_x]_j 
	\geq v_j + \frac{\Delta x_j}{2} [v_x]_j.
\end{gather}\end{subequations}
Next for \cref{eqn:reconWH_selfmon} we use
\begin{subequations}\begin{gather}
	\pdv{v_{j+1/2}^-}{v_j} = 1 + \frac{w_j^\downarrow}{w_{j+1/2}^-} \frac{\Delta x_j}{2} \pdv{[v_x]_j}{v_j}
	\qquad
	\text{and}
	\qquad
	\pdv{v_{j+1/2}^{v-}}{v_j} = 1 + \frac{\Delta x_j}{2} \pdv{[v_x]_j}{v_j},
\shortintertext{if $\pdv*{[v_x]_j}{v_j} \geq 0$ then}
	1
	\leq 1 + \frac{w_j^\downarrow}{w_{j+1/2}^-} \frac{\Delta x_j}{2} \pdv{[v_x]_j}{v_j}
	\leq  1 + \frac{\Delta x_j}{2} \pdv{[v_x]_j}{v_j},
\shortintertext{and if $\pdv*{[v_x]_j}{v_j} \leq 0$ then}
	1
	\geq 1 + \frac{w_j^\downarrow}{w_{j+1/2}^-} \frac{\Delta x_j}{2} \pdv{[v_x]_j}{v_j}
	\geq  1 + \frac{\Delta x_j}{2} \pdv{[v_x]_j}{v_j}.
\end{gather}\end{subequations}
Finally for \cref{eqn:reconWH_neighmon} we use
\begin{subequations}\begin{gather}
	\pdv{v_{j+1/2}^-}{v_{j+1}} = \frac{w_j^\downarrow}{w_{j+1/2}^-} \frac{\Delta x_j}{2} \pdv{[v_x]_j}{v_{j+1}}
	\qquad
	\text{and}
	\qquad
	\pdv{v_{j+1/2}^{v-}}{v_{j+1}} = \frac{\Delta x_j}{2} \pdv{[v_x]_j}{v_{j+1}},
\shortintertext{if $\pdv*{[v_x]_j}{v_{j+1}} \geq 0$ then}
	0
	\leq \frac{w_j^\downarrow}{w_{j+1/2}^-} \frac{\Delta x_j}{2} \pdv{[v_x]_j}{v_{j+1}}
	\leq  \frac{\Delta x_j}{2} \pdv{[v_x]_j}{v_{j+1}},
\shortintertext{and if $\pdv*{[v_x]_j}{v_j} \leq 0$ then}
	0
	\geq \frac{w_j^\downarrow}{w_{j+1/2}^-} \frac{\Delta x_j}{2} \pdv{[v_x]_j}{v_{j+1}}
	\geq  \frac{\Delta x_j}{2} \pdv{[v_x]_j}{v_{j+1}}.
\end{gather}\end{subequations}
\section{Reconstruction of concentration} \label{sec:supprecon_phi}

In this section we generalise the results of \cref{thm:results_phi,sec:supprecon_phi}. We define $\phi_j \eqdef [\phi h]_j / h_j$ and take
\begin{align}
[\phi _x]_j \eqdef \sigma^{\phi}(\phi_{j-1},\phi_{j},\phi_{j+1};\{x\}),
\end{align}
where $\sigma^{\phi}$ is a symmetric, TVD, and CDP slope limiter bounded as in \cref{thm:PWB_linearTVD} with the parameters $\alpha_{j+1/2}^{\pm}$ having least values $\alpha_{j+1/2}^{\phi\pm}$, $\alpha_j^{\phi\uparrow} \eqdef \max(\alpha_{j-1/2}^{\phi+},\alpha_{j+1/2}^{\phi-})$, and its derivative is bounded as in \cref{thm:PWB_selfmon} with $\beta_{j}^{\pm}$ having least values $\beta_{j}^{\phi\pm}$. We also define
\begin{equation}\begin{split}
	\phi_{j}^\downarrow &\eqdef \min \pbrk*{ \phi_j - \alpha_{j-1/2}^{\phi+}\Delta\phi_{j-1/2},\phi_j,\phi_j + \alpha_{j+1/2}^{\phi-}\Delta\phi_{j+1/2} },
	\\
	\phi_{j}^\uparrow &\eqdef \max \pbrk*{ \phi_j - \alpha_{j-1/2}^{\phi+}\Delta\phi_{j-1/2},\phi_j,\phi_j + \alpha_{j+1/2}^{\phi-}\Delta\phi_{j+1/2} }.
\end{split}\end{equation}

Exactly as was the case in \cref{sec:suppreconWH}, from $[\phi_x]_j$ we construct a gradient $[\Phi_x]_j$, then the values $[\Phi]_{j+1/2}^\pm$ by \cref{eqn:FV_recon}, and finally $\phi_{j+1/2}^\pm \eqdef [\Phi]_{j+1/2}^\pm/h_{j+1/2}^\pm$, these are the values to be used. Thus we may use the results of \cref{sec:reconWH}, and immediately set
\begin{equation}
	[\phi h_x]_j \eqdef [\phi_x]_j \ppar*{ h_j - \frac{\Delta x_j}{2} \abs*{[h_x]_j} } + \phi_j [h_x]_j.
\end{equation}
as our reconstruction, and thus
\begin{align}
	\phi_j^\downarrow &\leq \phi_{j-1/2}^+ \leq \phi_j^\uparrow,
	&
	\phi_j^\downarrow &\leq \phi_{j+1/2}^- \leq \phi_j^\uparrow.
\end{align}
However, unlike the case of \cref{sec:reconWH}, the reconstruction of $h$ is dependent on $\{\phi\}$ though $B_j$ \cref{eqn:B_fastfluid_concen}. Because $B_j$ is independent of $\phi_{j \pm 1}$ we immediately get neighbour-monotonicity
\begin{align}
	\pdv{\phi_{j-1/2}^+}{\phi_{j-1}} &\geq 0, &
	\pdv{\phi_{j+1/2}^-}{\phi_{j+1}} &\geq 0.
\end{align}
Self-monotonicity is more challenging because $B_j$ is a function of $\phi_j$. We consider $\phi_{j+1/2}^-$ ($\phi_{j-1/2}^+$ similar by symmetry), which is
\begin{align}
	\phi_{j+1/2}^- &= \phi_j + \frac{h_j - \frac{\Delta x_j}{2} \abs*{[h_x]_j}}{h_{j+1/2}^-} \frac{\Delta x_j}{2} [\phi_x]_j.
\end{align}
For $[h_x]_j \leq 0$ we have $\phi_{j+1/2}^- = \phi_j + [\phi_x]_j \frac*{\Delta x_j}{2}$, thus we have self-monotonicity $\pdv*{\phi_{j+1/2}^-}{\phi_{j}} \geq 0$. Also, when $\pdv*{[h_x]_j}{\phi_j} = 0$ then \cref{eqn:reconWH_selfmon} can be applied and, again, we have self-monotonicity. The only case left is $[h_x]_j \geq 0$, $1 \leq \xi_j \leq \xi_j^C$, and $\Delta b_j^\uparrow = B_j$, which we now consider. First,
\begin{gather*}
	\begin{aligned}
		\pdv{h_{j+1/2}^-}{\phi_j} &= G_j R_j \frac{h_{j}^\downarrow}{B_j} \pdv{B_j}{\phi_j},&
		\pdv{h_{j-1/2}^+}{\phi_j} &= -G_j R_j \frac{h_{j}^\downarrow}{B_j} \pdv{B_j}{\phi_j},
	\end{aligned}
\shortintertext{thus}
	\pdv{\phi_{j+1/2}^-}{\phi_j} 
	= 1 + \frac{\Delta x_j}{2} \frac{h_{j-1/2}^+}{h_{j+1/2}^-} \pdv{[\phi_x]_j}{\phi_j} - \frac{\Delta x_j}{2} \frac{[\phi_x]_j}{h_{j+1/2}^-} \cdot G_j R_j \frac{h_{j}^\downarrow}{B_j} \pdv{B_j}{\phi_j} \cdot \ppar*{ 1 + \frac{h_{j-1/2}^+}{h_{j+1/2}^-} }.
\end{gather*}
To construct a lower bound, we use that $\pdv*{[\phi_x]_j}{\phi_j} \geq -2 \beta_j^{\phi-}/\Delta x_j$ by \cref{thm:PWB_selfmon}, $\abs*{[\phi_x]_j} \leq 2 \alpha_j^{\phi\uparrow} \phi_j/\Delta x_j$ by \cref{thm:PWB_linearTVD}, $\abs*{R_j} \leq 1$ by \cref{thm:suppreconH_bound_R_minmod}, and define
\begin{gather}
	P_j = \sup \abs*{\frac{\phi_j}{B_j} \pdv{B_j}{\phi_j}} = \frac{1}{3}.
\shortintertext{Thus}
	\pdv{\phi_{j+1/2}^-}{\phi_j} \geq 1 - \beta_j^{\phi-} \frac{h_{j-1/2}^+}{h_{j+1/2}^-}  -  \alpha_j^{\phi\uparrow} G_j P_j \frac{h_{j}^\downarrow}{h_{j+1/2}^-} \ppar*{ 1 + \frac{h_{j-1/2}^+}{h_{j+1/2}^-} }.
\end{gather}
This expression contains a number of depth ratios, and because we consider $[h_x]_j \geq 0$ they all take a maximal value of $1$. Thus
\begin{align}
	\pdv{\phi_{j+1/2}^-}{\phi_j} &\geq 1 - \beta_j^{\phi-}  -  2 \alpha_j^{\phi\uparrow} G_j P_j .
\shortintertext{Therefore, if}
	G_j &\leq \frac{1 - \beta_j^{\phi-}}{2 \alpha_j^{\phi\uparrow} P_j } = \frac{3 \ppar*{1 - \beta_j^{\phi-}}}{2 \alpha_j^{\phi\uparrow} P_j }
\end{align}
then the reconstruction is self-monotone.
\section{Other results} \label{sec:suppother}

\begin{lemma} \label{thm:minmod_centdiff_condn}
	Suppose that $C = a L + b R$ where $a,b \geq 0$ and $L,C,R > 0$. 
	\begin{enumerate}[left=0mm]
	  \item If $C < \min(L,R)$, then $a+b < 1$ and one of $a,b$ is non-zero.
	  \item If $L < \min(C,R)$ then $a > 1$ or $b>0$ (or both).
	  \item If $R < \min(L,C)$ then $b > 1$ or $a>0$ (or both).
	\end{enumerate}
\end{lemma}
\begin{proof}\leavevmode
	\begin{enumerate}[left=0mm]
		\item If $L \leq R$ we have that $C < L$ which implies $(a-1)L < -bR$. If $a>1$ then $L < -bR/(a-1)$, not possible. If $a=1$ then $bR < 0$, not possible. If $a<1$ then $L > bR/(1-a)$, which implies that $b/(1-a)<1$ and $a+b < 1$. Applying the same logic to the case $L \geq R$ gives the same result, thus $a+b < 1$.
		\item $L < C = aL + bR$ thus $(1-a)L < bR$. If $a>1$ then $L > -bR(a-1)$, true. If $a=1$ then $0 < bR$, thus $b>0$. If $a<1$ then $L < bR(1-a)$, thus $b>0$.
		\item $R < C = aL + bR$ thus $(1-b)R < aL$. If $b>1$ then $R > -aL(b-1)$, true. If $b=1$ then $0 < aL$, thus $a>0$. If $b<1$ then $R < aL(1-b)$, thus $a>0$.
	\end{enumerate}
\end{proof}

\end{document}